\documentclass[runningheads]{llncs}
\usepackage{graphicx}
\usepackage{latexsym}
\usepackage{amssymb}
\usepackage{amsmath}
\usepackage{mathrsfs}
\usepackage{enumerate}
\begin{document}
\title{A unified relational semantics for BPL, IPL and OL -- axiomatization without disjunction
}
\titlerunning{A unified relational semantics for BPL, IPL and OL}
%
\author{Zhicheng Chen
}
\authorrunning{Z. Chen}
%
\institute{Department of Philosophy and Religious Studies, Peking University, Beijing, China 
}
\maketitle              
\begin{abstract}
In this paper, we propose a relational semantics of propositional language, which unifies the relational semantics of intuitionistic logic, Visser’s Basic Propositional Logic and orthologic. 
Working in language $\{\bot,\land,\neg\}$ and $\{\bot,\land,\to\}$ respectively, we axiomatize this basic logic as well as stronger ones corresponding to different combinations of frame conditions: reflexivity, symmetry, and transitivity. We also provide translations from these propositional logics into modal logics. 

\keywords{Unified semantics  \and Intuitionistic logic \and Basic Propositional Logic \and Orthologic.}
\end{abstract}
\section{Introduction}
Relational semantics is widely used in the study of non-classical propositional logics. Compared with its use in normal modal logics, a notable difference is that valuations of formulas do not range over the whole power set of the domain. For example, in intuitionistic propositional logic(\textbf{IPL}) and Visser's basic propositional logic(\textbf{BPL}), valuations are asked to be upward closed \cite{chagrov,visser}; in orthologic(\textbf{OL}), which is a weakened form of orthomodular quantum logic, valuations need to be bi-orthogonally closed \cite{goldblatt}.\par

On Pages 139-140 in \cite{dalla}, Dalla Chiara and Giuntini make an interesting observation that, in a Kripke model $\langle W,R,V\rangle$, if we require valuations to be:
\begin{align}
\Box RV(p) \subseteq V(p) \textit{, for each proposition letter } p, \tag{*} \label{*}
\end{align}
(where $R$ and $\Box$ are unary operators over $\wp(W)$, mapping $X\subseteq W$ to $RX=\{w\in W\mid \textit{there is an } x \in X \textit{ with }xRw\}$ and $\Box X=\{w\in W | R\{w\}\subseteq X\}$,)\\
then it is equivalent to the requirement of \textbf{IPL} when $R$ is reﬂexive and transitive, and equivalent to that of \textbf{OL} when $R$ is reflexive and symmetric. Throughout the paper, we shall use ``DG-models" to denote Kripke models satisfying \eqref{*}. 

Later, it is found (\cite{zhong.rg,wesley}) that the interpretation of connectives $\land,\lor$ and $\neg$ in relational semantics of \textbf{IPL} and \textbf{OL} can be unified in DG-models as follows:\\
In a DG-model $\langle W,R,V\rangle$, let
\begin{itemize}
    \item[] $\lVert\varphi\land\psi\rVert=\lVert\varphi\rVert\cap \rVert\psi\rVert$,
    \item[] $\lVert\neg\varphi\rVert = \Box-\lVert\varphi\rVert$,
    \item[] $\lVert\varphi\lor\psi\rVert=\Box R(\lVert\varphi\rVert\cup \lVert\psi\rVert)$,
\end{itemize}
where $\lVert\varphi\rVert$ is the set of possible worlds at which $\varphi$ holds, and $-X=W\setminus X$ for any $X\subseteq W$. 

Thus, there is a certain ``basic" logic lying beneath both \textbf{OL} and the $\{\land,\lor,\neg\}$-fragment of \textbf{IPL}. Let us call it ``DG logic". One can find the axiomatization for this logic in \cite{wesley}, and the axiomatization for the $\{\land,\neg\}$-fragment of this logic in \cite{zhong.rg} (--the valuation condition there is slightly different from \eqref{*}).

It should be mentioned that in \cite{wesley}, Holliday introduced a ``Fundamental Logic", which also lies beneath \textbf{IPL} and \textbf{OL}. It aims to capture the properties of connectives $\land,\lor$ and $\neg$ that are determined by their introduction and elimination rules in Fitch's natural deduction system. Fundamental Logic is stronger than DG logic, since it is sound and complete with respect to reflexive and ``pseudosymmetric" DG-models.

In this paper, we slightly strengthen the requirement \eqref{*} into:
\begin{align}
RV(p)\cap\Box RV(p) \subseteq V(p) \textit{, for each proposition letter } p, \tag{**} \label{**}
\end{align}
The effect is that valuations will be upward closed once $R$ is transitive. So we are able to bring \textbf{BPL} into the picture. 
We call Kripke models satisfying \eqref{**} as ``BIO-models".
And a new ``basic" logic emerges, which lies under \textbf{OL} and the $\{\land,\lor,\neg\}$-fragment of \textbf{BPL}, \textbf{IPL}, given the interpretation of connectives over BIO-models as follows:
\begin{itemize}
    \item[] $\lVert\varphi\land\psi\rVert=\lVert\varphi\rVert\cap \rVert\psi\rVert$,
    \item[] $\lVert\neg\varphi\rVert = \Box-\lVert\varphi\rVert$,
    \item[] $\lVert\varphi\lor\psi\rVert=\lVert\varphi\rVert\cup \lVert\psi\rVert\cup(R(\lVert\varphi\rVert\cup \lVert\psi\rVert)\cap\Box R(\lVert\varphi\rVert\cup \lVert\psi\rVert))$.
\end{itemize}
We call this logic ``BIO logic". 

The reader might have wondered why implication is left outside the framework.
The main reason is that there is no fixed definition for implication in orthomodular quantum logic, as well as \textbf{OL}. 

In \cite{dalla}, several implication candidates for orthomodular quantum logic and \textbf{OL} are discussed. 
One of the candidates is the ``Sasaki hook", whose interpretation in Kripkean semantics is as follows:
$$\lVert\varphi\hookrightarrow\psi\rVert=\Box (-\lVert\varphi\rVert\cup\Diamond (\lVert\varphi\rVert\cap\lVert\psi\rVert)),$$
where $\Diamond X=\{w\in W\mid\textit{ there is a }v\in X \textit{ with }wRv\}$.
If one wants to equip DG logic (or BIO logic) with an implication $\to$ to unify \textbf{IPL} and \textbf{OL} with Sasaki hook, then a natural way to interpret $\to$ over DG-models (or BIO-models) would be:
$$\lVert\varphi\to\psi\rVert=\Box (-\lVert\varphi\rVert\cup R(\lVert\varphi\rVert\cap\lVert\psi\rVert)).$$ 
In \cite{wesley.preconditional}, Holliday thoroughly studies the behavious of this kind of implication over DG-models.

Among other candidates mentioned in \cite{dalla} is the strict implication, which is exactly the kind defined in relational semantics of \textbf{IPL} and \textbf{BPL}. It is argued in \cite{Kawano} that in \textbf{OL}, strict implication has certain advantages over other candidates. 
Therefore, it seems desirable to consider BIO logic equipped with strict implication, which lies beneath \textbf{BPL}, \textbf{IPL} and \textbf{OL} with strict implication. 

In this paper, we will provide the axiomatization of the $\{\bot,\land,\neg\}$ fragment of BIO logic, and the $\{\bot,\land,\to\}$ fragment of BIO logic with strict implication--negation can be defined using $\to$ and $\bot$. We will also consider stronger logics arising from various combinations of three frame properties: reflexivity, symmetry, and transitivity.
However, the axiomatization of the full BIO logic (with disjunction included) turns out to be quite complicated, unlike the situation in \cite{wesley}. We defer this study to a follow-up article.


It should be mentioned that weak logics with strict implication have been well studied in the literature \cite{corsi,closer,bou,Ma Zhao}. A common feature of these studies is that, in addition to implication, the formal language contains local disjunction-``local" from the perspective of its interpretation in relational semantics. 
However, when local disjunction is absent, as in the unified semantics we are going to investigate, difficulties in axiomatization and completeness proof arise.

For example, suppose that we are working on a Gentzen-style system for strict implication with respect to the class of symmetric pointed models. With local disjunction, symmetry can be characterized by the following sequent: $$\{\varphi\}\vdash\psi\lor((\varphi\to\psi)\to\bot)$$

Note that there is only one formula, rather than a set of formulas on the right side. In fact, in this case we can establish for this logic a single succedent sequent calculus. 
If we do not have local disjunction, then at first sight we seem forced to use multisuccedent sequent  to express $\{\varphi\}\vdash\psi\lor((\varphi\to\psi)\to\bot)$ by
$\{\varphi\}\vdash\{\psi,((\varphi\to\psi)\to\bot)\}$. But is it still possible to have a single succedent sequent calculus for axiomatization? The same problem is faced in orthologic with strict implication, for which there are no single succedent sequent calculi in the literature, but only multisuccedent sequent calculi \cite{Kawano,Kawano.labeled}. 

The lack of local disjunction also leads to difficulties in proving completeness of logics with respect to relational semantics. 
With local disjunction, there’s seemingly a uniform definition of canonical model \cite{closer,bou}, by using sets of formulas each of which is consistent, deduction closed and “disjunction complete”: 
\[\textit{for any formula }\varphi,\psi\textit{, if }\varphi\lor\psi\in\Gamma\textit{, then }\varphi\in\Gamma\textit{ or }\psi\in\Gamma\] 
When local disjunction is absent, things become tricky. In some cases, we can get through this simply by removing the requirement of being ``disjunction complete". However, this does not seem to work in cases involving symmetry. 

We will solve these difficulties. More precisely, we will provide a sound and complete single succedent sequent calculus for orthologic with strict implication. And by generalizing the notion of being ``disjunction complete", we contribute to a uniform definition of the underlying set of the canonical model.

The rest of this paper is organized as follows. In Section 2 we provide the formal semantics and make some technical preparations. In Section 3, we axiomatize the $\{\bot,\land,\neg\}$ fragment of BIO logic as well as stronger logics induced by combinations of frame conditions: reflexivity, symmetry, and transitivity. In Section 4, we redo this task, but working in the language $\{\bot,\land,\to\}$. Finally in Section 5, translations from these propositional logics into modal logics are studied.

\section{Formal Semantics}

We assume that the reader is familiar with basic notions in Kripke semantics, such as the notion of frames, valuations, models, pointed models, and so on.


First, let's define some set operators that we will frequently use later: 

\begin{definition}
For any Kripke frame \( \mathfrak{F} = \langle W, R \rangle \), any \( X, Y \subseteq W \), define:
\begin{itemize}
    \item[] $R(X)=\{w\in W\mid \textit{there is an } x \in X \textit{ with }xRw\}$,
    \item[] $\Box_{\mathfrak{F}}(X) = \{ w \in W \mid R(\{w\}) \subseteq X \}$,
    \item[] $-_{\mathfrak{F}}(X) = W \setminus X$.
\end{itemize}

\end{definition}
(Note: sometimes for convenience, we will omit the subscript “\(\mathfrak{F}\)” or parentheses.)

It's easy to check that: 

\begin{lemma}
~
\begin{enumerate}[(1)]
    \item $R$ is distributive over $\cup$; $\Box_{\mathfrak{F}}$ is distributive over $\cap$.
    \item Both $R$ and $\Box_{\mathfrak{F}}$ are monotone.
    \item $RX\subseteq Y\Leftrightarrow X\subseteq \Box_{\mathfrak{F}}Y$.
    \item $\Box_{\mathfrak{F}}R\Box_{\mathfrak{F}}X=\Box_{\mathfrak{F}}X$; $R\Box_{\mathfrak{F}}RX=RX$.
\end{enumerate}
\end{lemma}

Let $PL=\{p_0,p_1,...\}$ be a countable set of propositional letters. 
The formulas we consider in this paper are those built from the logical constant $\bot$ and propositional letters, using connectives $\land, \neg, \to$. Denote by $Form$ the set of all formulas. For convenience, let us agree that the order of precedence of connectives is $\neg>\land>\,\to$.

\begin{definition}[Interpretation of connectives and formulas]
Given any Kripke frame \( \mathfrak{F} \), the connectives are interperted as certain set operators
in \( \mathfrak{F} \):
\begin{itemize}
    \item $\land$ is interpreted as $\cap$;
    \item $\neg$ is interpreted as an unary operator $\neg_{\mathfrak{F}}$ :
    for any $X\subseteq W$, 
    $\neg_{\mathfrak{F}} X = \Box_{\mathfrak{F}} -_{\mathfrak{F}}X$;
    \item $\to$ is interpreted as a binary operator $\to_{\mathfrak{F}}$ :
    for any $X,Y\subseteq W$, $X \to_{\mathfrak{F}} Y = \Box_{\mathfrak{F}} (-_{\mathfrak{F}}X \cup Y)$.
\end{itemize}

Then given any Kripke model \( \mathfrak{M}=(\mathfrak{F},V) \), we can define the truth set $\| \phi \|_{\mathfrak{M}}$ of formula $\phi$  recursively:
\begin{itemize}
    \item $\| \bot \|_{\mathfrak{M}} = \emptyset$;
    \item $\| p \|_{\mathfrak{M}} = V(p)$, for any $p \in PL$;
    \item $\| \alpha \land \beta \|_{\mathfrak{M}} = \| \alpha \|_{\mathfrak{M}} \cap \| \beta \|_{\mathfrak{M}}$, for any formula $\alpha,\beta$;
    \item $\|\neg \alpha \|_{\mathfrak{M}} = \neg_{\mathfrak{F}} \| \alpha \|_{\mathfrak{M}}$, for any formula $\alpha$;
    \item $\|\alpha \to \beta\|_{\mathfrak{M}} = \| \alpha \|_{\mathfrak{M}}\to_{\mathfrak{F}} \| \beta \|_{\mathfrak{M}}$, for any formula $\alpha,\beta$.
\end{itemize}

Then the satisfaction relation is defined by:
\[
\mathfrak{M}, w \models \phi \iff w \in \| \phi \|_{\mathfrak{M}},
\]
for any Kripke pointed model \( \mathcal{M}, w \) and formula \( \phi \).
\end{definition}

For the definition of BIO models, we introduce a key operator:

\begin{definition}
For any Kripke frame \( \mathfrak{F} = \langle W, R \rangle \), any \( X\subseteq W \), define
$$\text{BIO}_{\mathfrak{F}}(X) = X \cup (RX \cap \Box_{\mathfrak{F}}RX).$$
\end{definition}

\begin{lemma}
For any Kripke frame \( \mathfrak{F}= \langle W, R \rangle \), we have: 
    \begin{enumerate}[(1)]
        \item for any \( X\subseteq W \), $\text{BIO}_{\mathfrak{F}}(X) = (X \cup RX) \cap \Box_{\mathfrak{F}}RX$;
        \item $\text{BIO}_{\mathfrak{F}}$ is a closure operator. (i.e. $X\subseteq\text{BIO}_{\mathfrak{F}}(X)$, $X\subseteq Y\Rightarrow\text{BIO}_{\mathfrak{F}}(X)\subseteq\text{BIO}_{\mathfrak{F}}(Y)$, and $\text{BIO}_{\mathfrak{F}}(\text{BIO}_{\mathfrak{F}}(X))\subseteq\text{BIO}_{\mathfrak{F}}(X)$, for any \( X,Y\subseteq W \).) 
    \end{enumerate}
\end{lemma}

\begin{proof}
    \begin{enumerate}[(1)]
        \item It follows from $X\subseteq \Box_{\mathfrak{F}} RX$.
        \item We only prove $\text{BIO}_{\mathfrak{F}}(\text{BIO}_{\mathfrak{F}}(X))\subseteq\text{BIO}_{\mathfrak{F}}(X)$. It suffices to show $R(\text{BIO}_{\mathfrak{F}}(X))\subseteq RX$. From (1) we can infer that $BIO_{\mathfrak{F}}(X)\subseteq \Box_{\mathfrak{F}}RX$. So $R(\text{BIO}_{\mathfrak{F}}(X))\subseteq R\Box_{\mathfrak{F}}RX=RX$.
    \end{enumerate}
\end{proof}

\begin{definition}
    Let \( \mathfrak{F} = \langle W, R \rangle \) be a Kripke frame, and $X\subseteq W$. $X$ is called a fixed point of $\text{BIO}_{\mathfrak{F}}$, if $\text{BIO}_{\mathfrak{F}}(X)=X$. Denote by $FP(\text{BIO}_{\mathfrak{F}})$ the set of fixed points of $\text{BIO}_{\mathfrak{F}}$.
\end{definition}

\begin{lemma}\label{frame condition.fixpoint}
    Let \( \mathfrak{F} = \langle W, R \rangle \) be a Kripke frame.
    \begin{enumerate}[(1)]
        \item $R$ is reflexive $\ \Rightarrow\ $ $FP(\text{BIO}_{\mathfrak{F}})=\{X\subseteq W\mid\Box_{\mathfrak{F}} RX\subseteq X\}$.
        \item $R$ is symmetric $\ \Rightarrow\ $ $FP(\text{BIO}_{\mathfrak{F}})=\{X\subseteq W\mid\neg_{\mathfrak{F}}\neg_{\mathfrak{F}} X\subseteq \neg_{\mathfrak{F}} X \cup X\}$$=\{X\subseteq W\mid X\subseteq \neg_{\mathfrak{F}}\neg_{\mathfrak{F}} X \to_{\mathfrak{F}} X\}$.
        \item $R$ is transitive $\ \Rightarrow\ $ $FP(\text{BIO}_{\mathfrak{F}})=\{X\subseteq W\mid X\subseteq\Box_{\mathfrak{F}}X\}$.
        \item $R$ is symmetric and transitive $\ \Rightarrow\ $ $FP(\text{BIO}_{\mathfrak{F}})=\{X\subseteq W\mid W=X\cup\neg_{\mathfrak{F}}X\}$.
    \end{enumerate}
\end{lemma}

\begin{definition}[BIO model]
  A Kripke model $\mathfrak{M}=\langle \mathfrak{F},V\rangle$ is called a BIO model if for each $p\in PL$, $V(p)\in FP(\text{BIO}_{\mathfrak{F}})$. 
\end{definition}

\begin{lemma}\label{formula.fixedpoint}
Let $\mathfrak{M}=\langle \mathfrak{F},V\rangle$ be a BIO model. For each formula $\varphi$, $\Vert \varphi \Vert\in FP(\text{BIO}_{\mathfrak{F}})$.
\end{lemma}
\begin{proof}
Prove by induction on formulas, using the following facts: 
\begin{itemize}
    \item $\emptyset\in FP(\text{BIO}_{\mathfrak{F}})$; 
    \item $FP(\text{BIO}_{\mathfrak{F}})$ is closed under taking intersection;
    \item $\Box_{\mathfrak{F}}(X)\in FP(\text{BIO}_{\mathfrak{F}})$ for any $X\subseteq W$.
\end{itemize}
\end{proof}

\begin{definition}\label{classofpi}
For any $\spadesuit\in\{\mathbf{K},\mathbf{T}, \mathbf{B}, \mathbf{4}, \mathbf{B4}, \mathbf{S4}, \mathbf{TB}, \mathbf{S5}\}$, let $\mathcal{D}_\spadesuit$ be the class of Kripke models that is characterized by the modal logic system $\spadesuit$.  Define 
\[\mathcal{D}_{\spadesuit^p}=\mathcal{D}_{\spadesuit}\cap\mathcal{D}_{\mathbf{K}^p},\]
where $\mathcal{D}_{\mathbf{K}^p}=$ $\{\,\mathfrak{M}\;|\;\mathfrak{M}$ is a BIO model $\}$.

Moreover, let
$\mathcal{D}_{\mathbf{I}}=\mathcal{D}_{\mathbf{S4}^p}$, $\mathcal{D}_{\mathbf{O}}=\mathcal{D}_{\mathbf{TB}^p}$, $\mathcal{D}_{\mathbf{C}}=\mathcal{D}_{\mathbf{S5}^p}$. ($\mathbf{V}$, $\mathbf{I}$, $\mathbf{O}$ and $\mathbf{C}$ stand for Visser's \textbf{BPL}, \textbf{IPL}, \textbf{OL} and classical propositional logic, respectively.)

For technical purpose, we introduce two additional classes:
\begin{itemize}
    \item[]  $\mathcal{D}_{\mathbf{K}^p -}=$ $\{\,\mathfrak{M}\;|\;\mathfrak{M}=(\mathfrak{F},V)=(W,R,V)$ is a Kripke model such that for any $p\in PL$, $RV(p) \cap\Box_{\mathfrak{F}}\emptyset \subseteq V(p)$ $\}$,
    \item[]  $\mathcal{D}_{\mathbf{T}^p -}=$ $\{\,\mathfrak{M}\;|\;\mathfrak{M}$ is a reflexive Kripke model$\}$.

\end{itemize}

\end{definition}

It is easy to verify that $\mathcal{D}_{\mathbf{K}^p}\subseteq\mathcal{D}_{\mathbf{K}^p-}$ and $\mathcal{D}_{\mathbf{T}^p}\subseteq\mathcal{D}_{\mathbf{T}^p-}$. And similar to lemma \ref{formula.fixedpoint}, we can prove by induction that:
\begin{lemma}\label{formula.fixedpoint.K-}
For any $\mathfrak{M}=(\mathfrak{F},V)=(W,R,V)\in\mathcal{D}_{\mathbf{K}^p -}$ and any formula $\alpha$, $R\Vert\alpha\Vert_{\mathfrak{M}} \cap\Box_{\mathfrak{F}}\emptyset \subseteq \Vert\alpha\Vert_{\mathfrak{M}}$.
\end{lemma}

\begin{definition}[Semantic consequence]
    For any set \( \Gamma \) of formulas, any formula \( \phi \), for each $\clubsuit\in\{\mathbf{K}^p, \mathbf{T}^p, \mathbf{B}^p, \mathbf{V}, \mathbf{B4}^p, \mathbf{I}, \mathbf{O}, \mathbf{C}, \mathbf{K}^p-, \mathbf{T}^p-\}$,  define:
\[
\Gamma \vDash_\clubsuit \phi \iff \text{for any } \mathfrak{M} \in \mathcal{D}_\clubsuit, \text{ for any } w \in \mathfrak{M}, (\mathfrak{M}, w \vDash \Gamma \Rightarrow \mathfrak{M}, w \vDash \phi),\]
where $\mathfrak{M}, w \vDash \Gamma$ means $\mathfrak{M}, w \vDash \alpha$ for each $\alpha\in\Gamma$.
\end{definition}

It follows from $\mathcal{D}_{\mathbf{K}^p}\subseteq\mathcal{D}_{\mathbf{K}^p-}$ and $\mathcal{D}_{\mathbf{T}^p}\subseteq\mathcal{D}_{\mathbf{T}^p-}$ that, $\,\vDash_{\mathbf{K}^p-}\subseteq\,\vDash_{\mathbf{K}^p}$ and $\,\vDash_{\mathbf{T}^p-}\subseteq\,\vDash_{\mathbf{T}^p}$.

Somewhat surprisingly, we will show that the reverse inclusion is also true, so that the semantic consequence relation $\,\vDash_{\mathbf{K}^p}$ actually equals $\,\vDash_{\mathbf{K}^p-}$, and $\,\vDash_{\mathbf{T}^p}$ equals $\,\vDash_{\mathbf{T}^p-}$.
This phenomenon is essentially due to the limited expressive power of our language.

First, let us review a method for model transformation called ``unravelling".

\begin{definition}
~
\begin{itemize}
  \item Let $\mathfrak{M},w=(W,R,V,w)$ be a pointed Kripke model. The unravelling of $\mathfrak{M}$ from $w$ is the model $Unr_w(\mathfrak{M})=(W',R',V')$, where
    \begin{enumerate}[(i)]
        \item $W'=\{(s_0,...,s_n)\ |\ n\in \mathbb{N}, s_0,...,s_n\in W$ satisfying that $s_0=w$ and $s_i Rs_{i+1}$ for each $i=0,...,n-1\}$;
        \item $(s_0,...,s_n)R'(t_0,...,t_m)$  $\iff m=n+1$ and $(s_0,...,s_n)=(t_0,...,t_n)$;
        \item  $V'(p)=\{(s_0,...,s_n)\in W'\ |\ s_n\in V(p)\}$ for each $p\in PL$.
    \end{enumerate}

    \item The reflexive unravelling of $\mathfrak{M}$ from $w$ is the model $Unr^{re}_w(\mathfrak{M})=(W',$ $R^*,V')$, where $W'$ and $V'$ are defined above, and $R^*$ is the reflexive closure of $R'$.
\end{itemize}
\end{definition}

\begin{lemma}\label{basic unravelling}
Let $\mathfrak{M},w$ be a pointed model. 
\begin{enumerate}[(1)]
    \item For any $\varphi\in Form$ and $(s_0,...,s_n)\in W'$, 
    $$Unr_w(\mathfrak{M}),(s_0,...,s_n)\vDash\varphi \iff \mathfrak{M},s_n\vDash\varphi.$$
    
    \item If $\mathfrak{M}\in\mathcal{D}_{\mathbf{K}^p-}$, then $Unr_w(\mathfrak{M})\in\mathcal{D}_{\mathbf{K}^p}$.

    \item If $\mathfrak{M}\in\mathcal{D}_{\mathbf{T}^p-}$, then for any $\varphi\in Form$ and $(s_0,...,s_n)\in W'$,
    $$Unr^{re}_w(\mathfrak{M}),(s_0,...,s_n)\vDash\varphi \iff \mathfrak{M},s_n\vDash\varphi.$$

    \item If $\mathfrak{M}\in\mathcal{D}_{\mathbf{T}^p-}$, then  $Unr^{re}_w(\mathfrak{M})\in\mathcal{D}_{\mathbf{T}^p}$.
\end{enumerate}
\end{lemma}

\begin{proof}~
\begin{enumerate}
    \item Let $f: W'\to W$ be defined by $f(s_0,...,s_n)=s_n$ for any $(s_0,...,s_n)\in W'$. It is known that $f$ is a bounded morphism from $Unr_w(\mathfrak{M})$ to $\mathfrak{M}$ (Lemma 4.52 in \cite{Blackburn}), and the modal satisfaction is invariant under bounded morphisms (Proposition 2.14 in \cite{Blackburn}). Notice that our propositional formulas can be expressed in modal language - for example, $\neg\alpha$ is expressed as $\Box\neg^c\alpha$ and $\alpha\to\beta$ is expressed as $\Box(\alpha\to^c\beta)$, where $\neg^c$ and $\to^c$ respectively denote the classical negation and classical implication in modal language.    
    As a consequence, for each of our propositional formulas, its satisfaction is also invariant under $f$.
    
    \item Assume $\mathfrak{M}=(\mathfrak{F},V)=(W,R,V)\in\mathcal{D}_{\mathbf{K}^p-}$. By the definition of $\mathcal{D}_{\mathbf{K}^p-}$, for any $\alpha\in Form_{\neg}$, $R\Vert\alpha\Vert_{\mathfrak{M}} \cap\Box_{\mathfrak{F}}\emptyset \subseteq \Vert\alpha\Vert_{\mathfrak{M}}$. Let $Unr_w(\mathfrak{M})=(\mathfrak{F}',V')=(W',R',V')$. We want to show that for all $p\in PL$, $R'V'(p) \cap\Box_{\mathfrak{F}'}R'V'(p) \subseteq V'(p)$. For any $p\in PL$ and $(s_0,...,s_n)\in W'$,
    
    \indent\quad $(s_0,...,s_n)\in R'V'(p) \cap\Box_{\mathfrak{F}'}R'V'(p)$ 
    
    $\Leftrightarrow$ there is a $R'$-predecessor of $(s_0,...,s_n)$ in $V'(p)$, and for any $R'$-successor of $(s_0,...,s_n)$, it has a $R'$-predecessor in $V'(p)$  \hfill (definition of the set operators $R'$ and $\Box_{\mathfrak{F}'}$)
    
    $\Leftrightarrow$ $s_{n-1}\in V(p)$, and for any $R'$-successor of $(s_0,...,s_n)$, it has a $R'$-predecessor in $V'(p)$  \hfill (the only $R'$-predecessor of $(s_0,...,s_n)$ is $(s_0,...,s_{n-1})$;\\ definition of $V'$)
    
    $\Leftrightarrow$ $s_{n-1}\in V(p)$, and, either $(s_0,...,s_n)$ has no $R'$-successor or $(s_0,...,s_n)\in V'(p)$  \hfill \quad(for each $R'$-successor of $(s_0,...,s_n)$, if any, it has only one $R'$-predecessor, i.e. $(s_0,...,s_n)$)
    
    $\Rightarrow$ $(s_0,...,s_n)\in V'(p)$.
    
    For the last step: Assume that $s_{n-1}\in V(p)$ and $(s_0,...,s_n)$ have no $R'$-successor. Then $s_n$ has no $R$-successor in $\mathfrak{M}$, according to the definition of $W'$ and $R'$. Then $s_n\in \Box_{\mathfrak{F}}\emptyset$ by the definition of $\Box_{\mathfrak{F}'}$. Since $s_{n-1}\in V(p)$ and $s_{n-1}Rs_n$ and $RV(p) \cap\Box_{\mathfrak{F}}\emptyset \subseteq V(p)$, we have $s_n\in V(p)$. So $(s_0,...,s_n)\in V'(p)$ by the definition of $V'$.

    \item Similar to (1). Given that $\mathfrak{M}$ is a reflexive model, $f$ can be proved to be a bounded morphism from $Unr^{re}_w(\mathfrak{M})$ to $\mathfrak{M}$.

    \item Assume that $\mathfrak{M}=(\mathfrak{F},V)=(W,R,V)\in\mathcal{D}_{\mathbf{T}^p-}$. Thus $R$ is reflexive. Let $Unr^{re}_w(\mathfrak{M})=(\mathfrak{F}^*,V')=(W',R^*,V')$. By definition $R^*$ is reflexive. Now we show that $Unr^{re}_w(\mathfrak{M})$ is a BIO-model. Since $R^*$ is reflexive, combining lemma \ref{frame condition.fixpoint}, it suffices to show that for all $p\in PL$, $\Box_{\mathfrak{F}^*}R^*V'(p) \subseteq V'(p)$.
    For any $p\in PL$ and $(s_0,...,s_n)\in W'$,
    
    \indent\quad $(s_0,...,s_n)\in  \Box_{\mathfrak{F}^*}R^*V'(p)$ 
    
    $\Leftrightarrow$ for any $R^*$-successor of $(s_0,...,s_n)$, it has a $R^*$-predecessor in $V'(p)$  \hfill (definition of the set operators $R^*$ and $\Box_{\mathfrak{F}^*}$)
    
    $\Rightarrow$  $(s_0,...,s_n,s_n)$ has a $R^*$-predecessor in $V'(p)$  \hfill (since $s_nRs_n$ ($R$ is reflexive), $(s_0,...,s_n,s_n)\in W'$ and $(s_0,...,s_n)R^*(s_0,...,s_n,s_n)$ )
    
    $\Leftrightarrow$ $(s_0,...,s_n,s_n)\in V'(p)$ or $(s_0,...,s_n)\in V'(p)$   \hfill ($(s_0,...,s_n,s_n)$ has only two  $R^*$-predecessors, itself and $(s_0,...,s_n)$)
    
    $\Leftrightarrow$ $s_n\in V(p)$  \hfill \quad(definition of $V'$).
\end{enumerate}
\end{proof}

Now we can show the desired equivalence.
\begin{lemma}\label{vDash_K-=vDash_K, vDash_T-=vDash_T}
 For any set \( \Gamma \) of formulas, any formula \( \phi \), 
 \begin{align*}
     \Gamma\,\vDash_{\mathbf{K}^p-}\phi \iff\Gamma\,\vDash_{\mathbf{K}^p}\phi\,;\\
     \Gamma\,\vDash_{\mathbf{T}^p-}\phi \iff\Gamma\,\vDash_{\mathbf{T}^p}\phi\,.
 \end{align*}
\end{lemma}
\begin{proof}
Let $\clubsuit\in\{\mathbf{K}^p,\mathbf{T}^p\}$ and fix $\Gamma$ and $\phi$. We want to show that $\Gamma\,\vDash_{\clubsuit-}\phi \Leftrightarrow\Gamma\,\vDash_{\clubsuit}\phi$. The non-trivial direction is from right to left, and we prove the contrapositive.

Assume that $\Gamma\nvDash_{\clubsuit-}\varphi$. Then by definition, there exist $\mathfrak{M}\in \mathcal{D}_{\clubsuit-}$ and $w\in\mathfrak{M}$ such that $\mathfrak{M},w\vDash\Gamma$ and $\mathfrak{M},w\nvDash\varphi$. Now let $\mathfrak{N}$ be: if $\clubsuit=\mathbf{K}^p$, the unravelling of $\mathfrak{M}$ from $w$; else if $\clubsuit=\mathbf{T}^p$, the reflexive unravelling of $\mathfrak{M}$ from $w$. By lemma \ref{basic unravelling}, we obtain $\mathfrak{N},w\vDash\Gamma$, $\mathfrak{N},w\nvDash\varphi$, and $\mathfrak{N}\in\mathcal{D}_{\clubsuit}$. So $\Gamma\nvDash_{\clubsuit}\varphi$.
\end{proof}

\section{Axiomatization in $\{\bot,\land,\neg\}$-language}

In this section we axiomatize $\vDash_{\mathbf{K}^p}$, $\vDash_{\mathbf{V}}$, $\vDash_{\mathbf{B}^p}$, and $\vDash_{\mathbf{B4}^p}$ in the $\{\bot,\land,\neg\}$-language. For the axiomatization of $\vDash_{\mathbf{T}^p}$, $\vDash_{\mathbf{I}}$, $\vDash_{\mathbf{O}}$, and $\vDash_{\mathbf{C}}$ in this language, see \cite{zhong}. For each system, we will provide a single succedent sequent calculus, which in our sense is a certain formal theory of the derivability relation (between sets of formulas and individual formulas).

Let us make the convention within this section that $T$ stands for $\neg\bot$.
Moreover, let $Form_\neg$ collect the formulas without the occurrence of $\to$, which will be the only formulas that we consider in this section.

The following concept is about how we express the classical disjunction in our restricted language.

\begin{definition}[``$\phi \vdash (\alpha_1, \ldots, \alpha_n)$"]
    For any $\,\vdash\,\subseteq \wp(Form_\neg) \times Form_\neg$, for any \( 2\leq n \in \mathbb{N} \), and for any \( \phi, \alpha_1, \ldots, \alpha_n\in Form_\neg \), we denote ``$\phi \vdash (\alpha_1, \ldots, \alpha_n)$" if and only if :
    \begin{enumerate}[(i)]
        \item for any $\Gamma\subseteq Form_\neg$ and $\psi\in Form_\neg$, if $\Gamma \cup \{ \alpha_i \} \vdash \psi$ for any $i \in \{ 1, \ldots, n \}$, then $\Gamma \cup \{ \phi \} \vdash \psi$;
        \item for any $\psi\in Form_\neg$, $\neg(\psi \land \alpha_1) \land \ldots \land \neg(\psi \land \alpha_n) \vdash \neg(\psi \land \phi)$.
    \end{enumerate}
\end{definition}

The idea of ``$\phi \vdash (\alpha_1, \ldots, \alpha_n)$" is that we aim to express ``$\phi \vdash \alpha_1 \lor^c \ldots \lor^c \alpha_n$", where $\lor^c$ represents the classical disjunction. However, given the current limited language, the most we can articulate are the two conditions above: one being the ``prove-by-cases" rule, and the other concerning the interplay between this classical disjunction and the negation in our language.

\begin{definition}[$\,\vDash_{\mathfrak{M},w}$]
    Let $\mathfrak{M},w$ be an arbitrary pointed Kripke model. Define $\,\vDash_{\mathfrak{M},w} \subseteq \wp(Form_\neg) \times Form_\neg$ to be: for any $\Gamma\subseteq Form_\neg$ and $\phi\in Form_\neg$,
    \[\Gamma\vDash_{\mathfrak{M},w}\phi\; \iff\; \mathfrak{M},w\vDash\Gamma \textit{ implies } \mathfrak{M},w\vDash\phi.\]
\end{definition}

\begin{lemma}\label{pseudo disjunction.soundness}
    For any Kripke model $\mathfrak{M}$, for any \( 2\leq n \in \mathbb{N} \), and for any \( \phi, \alpha_1, \ldots, \alpha_n\in Form_\neg \), if $\Vert\phi\Vert_{\mathfrak{M}} \subseteq \Vert \alpha_1\Vert_{\mathfrak{M}} \cup ... \cup\Vert\alpha_n \Vert_{\mathfrak{M}}$, then for any $w\in \mathfrak{M}$, $\phi \,\vDash_{\mathfrak{M},w} (\alpha_1, \ldots, \alpha_n)$. 
\end{lemma}

\begin{proof}
    By definition, $\phi \,\vDash_{\mathfrak{M},w} (\alpha_1, \ldots, \alpha_n)$ if and only if :
    \begin{enumerate}[(i)]
        \item for any $\Gamma\subseteq Form_\neg$ and $\psi\in Form_\neg$, if $\Gamma \cup \{ \alpha_i \} \,\vDash_{\mathfrak{M},w} \psi$ for any $i \in \{ 1, \ldots, n \}$, then $\Gamma \cup \{ \phi \} \,\vDash_{\mathfrak{M},w} \psi$;
        \item for any $\psi\in Form_\neg$, $\neg(\psi \land \alpha_1) \land \ldots \land \neg(\psi \land \alpha_n) \,\vDash_{\mathfrak{M},w} \neg(\psi \land \phi)$.
    \end{enumerate}
    These can be readily validated via $\Vert\phi\Vert_{\mathfrak{M}} \subseteq \Vert \alpha_1\Vert_{\mathfrak{M}} \cup ... \cup\Vert\alpha_n \Vert_{\mathfrak{M}}$ and the definition of $\,\vDash_{\mathfrak{M},w}$.
\end{proof}

\subsection{Sequent calculi}
We outline several sequent rules that apply later in the axiomatizations. Instead of detailing the rules by themselves, we describe what it means for $\ \vdash\ \subseteq \wp(Form_\neg) \times Form_\neg$ to satisfy each rule, omitting initial universal quantifiers for simplicity.

\begin{itemize}
    \item[] (A)\quad \( \Gamma \cup \{ \phi \} \vdash \phi \);
    \item[] (Cut)\quad \( \Gamma \cup \{ \psi \} \vdash \phi \text{ and } \Delta \vdash \psi \Rightarrow \Gamma \cup \Delta \vdash \phi \);
    \item[] (\(\land\)I)\quad \( \{ \phi, \psi \} \vdash \phi \land \psi \);
    \item[] (\(\land\)E)\quad \( \phi \land \psi \vdash \phi \) and \( \phi \land \psi \vdash \psi  \);
    \item[] (\(\bot\))\quad \( \bot \vdash \phi \);
    \item[] (\(\neg \bot\))\quad  \( \vdash T\);
    \item[] (\(\neg\) antitone)\quad \( \alpha \vdash \beta \Rightarrow \neg \beta \vdash \neg \alpha \);
    ~\\
    \item[] (Prop${}_{\mathbf{K}^p-}$)\quad \( \alpha \land \neg(\psi \land \alpha) \vdash \neg(\psi \land \neg T) \);
    \item[] ($\neg\neg$I) \quad  \( \alpha \vdash \neg \neg \alpha \);
    \item[] (Prop${}_{\mathbf{B}^p}$)\quad $\neg \neg \alpha\vdash (\alpha, \neg \alpha)$ ;
    \item[] (Prop${}_{\mathbf{V}}$)\quad \( \alpha \land \neg(\psi \land \alpha) \vdash \neg \psi \);
    \item[] (Prop${}_{\mathbf{B4}^p}$)\quad \( T\vdash (\alpha, \neg \alpha) \) .
\end{itemize}

Rules from (A) to (\(\neg\) antitone) are called basic rules. (\(\neg \bot\)) and (\(\neg\) antitone) are called basic negation rules. 

Furthermore, here are some important derived rules or properties:
\begin{itemize}
    \item[] (Mon)\quad \( \Gamma \vdash \phi \text{ and } \Gamma\subseteq\Delta  \Rightarrow \Delta \vdash \phi \).
    \item[] (Com)\quad If $\Gamma\vdash\varphi$, then there is a finite $\Gamma'\subseteq\Gamma$ such that $\Gamma'\vdash\varphi$ . 
\end{itemize}

It's easy to check that if $\vdash$ satisfies (A) and (Cut), then $\vdash$ also satisfies (Mon).

\begin{definition}
  Let \( \vdash_{\mathbf{K}^p} \) (and $\vdash_{\mathbf{V}}$, $\vdash_{\mathbf{B}^p}$,  $\vdash_{\mathbf{B4}^p}$, respectively) be the smallest $\ \vdash\ \subseteq \wp(Form_\neg) \times Form_\neg$ that satisfies, in addition to the basic rules, (Prop${}_{\mathbf{K}^p-}$) (and (Prop${}_{\mathbf{V}}$), 
($\neg\neg$I)+(Prop${}_{\mathbf{B}^p}$), ($\neg\neg$I)+(Prop${}_{\mathbf{B4}^p}$), respectively).
\end{definition}

\begin{lemma}
    $\vdash_{\mathbf{K}^p}$, $\vdash_{\mathbf{V}}$, $\vdash_{\mathbf{B}^p}$, and $\vdash_{\mathbf{B4}^p}$ all enjoy (Com).
\end{lemma}

\begin{proof}
    Take $\vdash_{\textbf{K}^p}$ as example. Let $S=\{(\Gamma,\varphi)\ |\ $there is a finite $\Gamma'\subseteq\Gamma$ with $\Gamma'\vdash_{\textbf{K}^p}\varphi\}$. Then the target is equivalent to $\ \vdash_{\textbf{K}^p}\subseteq S$.
    It suffices to show that for each rule in the definition of $\vdash_{\textbf{K}^p}$, if all sequents in the premises belong to $S$, then the sequent in the conclusion is also in $S$. The verification is straightforward and details are omitted.
\end{proof}

\subsection{Soundness}

\begin{definition}[Correctness at a pointed model]
    A sequent rule $(X)$ is correct at a pointed Kripke model $\mathfrak{M},w$, if and only if $\ \vDash_{\mathfrak{M},w}$ satisfies $(X)$.
    
   
\end{definition}

It is easy to show that :
\begin{lemma}
    The basic rules are correct at any pointed Kripke model. And ($\neg\neg$I) is correct at any pointed symmetric model.
\end{lemma}

\begin{lemma}
For each $\clubsuit\in\{\mathbf{K}^p-, \mathbf{V}, \mathbf{B}^p, \mathbf{B4}^p\}$, for any $\mathfrak{M}\in \mathcal{D}_{\clubsuit}$ and $w\in\mathfrak{M}$, (Prop${}_{\clubsuit}$) is correct at $\mathfrak{M},w$.

\end{lemma}

\begin{proof}
In this proof, we omit the subscript "$\mathfrak{M}$" of $\Vert\cdot\Vert$. 
\begin{enumerate}[(1)]
    \item For any model $\mathfrak{M}=(\mathfrak{F},V)=(W,R,V)$,\\
    \indent\quad $\mathfrak{M}\in \mathcal{D}_{\mathbf{K}^p-}$ \\
    $\Leftrightarrow$ for any $\alpha\in Form_{\neg}$,  $R\Vert\alpha\Vert \cap\Box_{\mathfrak{F}}\emptyset \subseteq \Vert\alpha\Vert$  \hfill (definition of $\mathcal{D}_{\mathbf{K}^p-}$, lemma \ref{formula.fixedpoint.K-})\\
    $\Leftrightarrow$ for any $\alpha\in Form_{\neg}$, $\Vert \alpha \Vert \subseteq \Box_{\mathfrak{F}}(-_{\mathfrak{F}}\Vert \neg T\Vert\cup \Vert\alpha \Vert)$ \hfill ($RX\subseteq Y\Leftrightarrow X\subseteq \Box_{\mathfrak{F}} Y$, $\Vert\neg T\Vert=\Box_{\mathfrak{F}}\emptyset$)\\
    $\Rightarrow$ for any $\alpha,\psi\in Form_{\neg}$, $\Vert \alpha \Vert\cap\Box_{\mathfrak{F}}(-_{\mathfrak{F}}\Vert \psi\Vert\cup -_{\mathfrak{F}}\Vert \alpha\Vert) \subseteq \Box_{\mathfrak{F}}(-_{\mathfrak{F}}\Vert \psi\Vert\cup-_{\mathfrak{F}}\Vert \neg T\Vert)$  \hfill \quad($\Box_{\mathfrak{F}}$ is distributive over $\cap$)\\
    $\Leftrightarrow$ for any $\alpha,\psi\in Form_{\neg}$, $\Vert \alpha\land\neg(\psi\land\alpha) \Vert \subseteq \Vert\neg(\psi\land\neg T)\Vert$  \hfill \quad(definition of $\Vert\cdot\Vert$)\\
    $\Leftrightarrow$ for any $\alpha,\psi\in Form_{\neg}$, for any $w\in\mathfrak{M}$, $ \alpha\land\neg(\psi\land\alpha)  \,\vDash_{\mathfrak{M},w} \neg(\psi\land\neg T)$  \hfill \quad(definition of the satisfaction relation and $\,\vDash_{\mathfrak{M},w}$)\\
    $\Leftrightarrow$ for any $w\in\mathfrak{M}$, (Prop${}_{\mathbf{K}^p-}$) is correct at $\mathfrak{M},w$  \hfill \quad(definition of correctness).

    \item For any model $\mathfrak{M}=(\mathfrak{F},V)=(W,R,V)$,\\
    \indent\quad $\mathfrak{M}\in \mathcal{D}_{\mathbf{V}}$ \\
    $\Leftrightarrow$ for any $\alpha\in Form_{\neg}$, $\Vert \alpha \Vert \subseteq \Box_{\mathfrak{F}}\Vert\alpha \Vert$ \hfill (definition of $\mathcal{D}_{\mathbf{V}}$, lemma \ref{formula.fixedpoint}, lemma \ref{frame condition.fixpoint})\\
    $\Rightarrow$ for any $\alpha,\psi\in Form_{\neg}$, $\Vert \alpha \Vert\cap\Box_{\mathfrak{F}}(-_{\mathfrak{F}}\Vert \psi\Vert\cup -_{\mathfrak{F}}\Vert \alpha\Vert) \subseteq \Box_{\mathfrak{F}}-_{\mathfrak{F}}\Vert \psi\Vert$  \hfill \quad($\Box_{\mathfrak{F}}$ is distributive over $\cap$)\\
    $\Leftrightarrow$ for any $\alpha,\psi\in Form_{\neg}$, $\Vert \alpha\land\neg(\psi\land\alpha) \Vert \subseteq \Vert\neg\psi\Vert$  \hfill \quad(definition of $\Vert\cdot\Vert$)\\
    $\Leftrightarrow$ for any $\alpha,\psi\in Form_{\neg}$, for any $w\in\mathfrak{M}$, $ \alpha\land\neg(\psi\land\alpha)  \,\vDash_{\mathfrak{M},w} \neg\psi$  \hfill \quad(definition of the satisfaction relation and $\,\vDash_{\mathfrak{M},w}$)\\
    $\Leftrightarrow$ for any $w\in\mathfrak{M}$, (Prop${}_{\mathbf{V}}$) is correct at $\mathfrak{M},w$  \hfill \quad(definition of correctness).
    
    \item For any model $\mathfrak{M}=(\mathfrak{F},V)=(W,R,V)$,\\
    \indent\quad $\mathfrak{M}\in \mathcal{D}_{\mathbf{B}^p}$ \\
    $\Leftrightarrow$ for any $\alpha\in Form_{\neg}$, $\neg_{\mathfrak{F}}\neg_{\mathfrak{F}}\Vert\alpha\Vert \subseteq \neg_{\mathfrak{F}}\Vert \alpha\Vert\cup \Vert\alpha \Vert$ \hfill (definition of $\mathcal{D}_{\mathbf{B}^p}$, lemma \ref{formula.fixedpoint}, lemma \ref{frame condition.fixpoint})\\
    $\Leftrightarrow$ for any $\alpha\in Form_{\neg}$, $\Vert\neg\neg\alpha\Vert \subseteq \Vert \neg\alpha\Vert\cup \Vert\alpha \Vert$  \hfill \quad(definition of $\Vert\cdot\Vert$)\\
    $\Rightarrow$ for any $\alpha\in Form_{\neg}$, for any $w\in\mathfrak{M}$, $\neg\neg\alpha \,\vDash_{\mathfrak{M},w} (\alpha, \neg\alpha)$  \hfill \quad(lemma \ref{pseudo disjunction.soundness})\\
    $\Leftrightarrow$ for any $w\in\mathfrak{M}$, (Prop${}_{\mathbf{B}^p}$) is correct at $\mathfrak{M},w$  \hfill \quad(definition of correctness).

    \item For any model $\mathfrak{M}=(\mathfrak{F},V)=(W,R,V)$,\\
    \indent\quad $\mathfrak{M}\in \mathcal{D}_{\mathbf{B4}^p}$ \\
    $\Leftrightarrow$ for any $\alpha\in Form_{\neg}$, $W \subseteq \neg_{\mathfrak{F}}\Vert \alpha\Vert\cup \Vert\alpha \Vert$ \hfill (definition of $\mathcal{D}_{\mathbf{B}^p}$, lemma \ref{formula.fixedpoint}, lemma \ref{frame condition.fixpoint})\\
    $\Leftrightarrow$ for any $\alpha\in Form_{\neg}$, $\Vert T \Vert \subseteq \Vert \neg\alpha\Vert\cup \Vert\alpha \Vert$  \hfill \quad(definition of $\Vert\cdot\Vert$, $\Vert T \Vert=W$)\\
    $\Rightarrow$ for any $\alpha\in Form_{\neg}$, for any $w\in\mathfrak{M}$, $T\,\vDash_{\mathfrak{M},w} (\alpha, \neg\alpha)$  \hfill \quad(lemma \ref{pseudo disjunction.soundness})\\
    $\Leftrightarrow$ for any $w\in\mathfrak{M}$, (Prop${}_{\mathbf{B4}^p}$) is correct at $\mathfrak{M},w$  \hfill \quad(definition of correctness).
   
\end{enumerate}
\end{proof}

\begin{theorem}[Soundness]
For each $\clubsuit\in\{\mathbf{K}^p, \mathbf{V}, \mathbf{B}^p, \mathbf{B4}^p\}$, for any $\Gamma\subseteq Form_\neg$ and $\phi\in Form_\neg$,
    \[\Gamma\,\vdash_\clubsuit\phi \Rightarrow\Gamma\,\vDash_\clubsuit\phi.\]
\end{theorem}
\begin{proof}
    Fix $\clubsuit\in\{\mathbf{K}^p, \mathbf{V}, \mathbf{B}^p, \mathbf{B4}^p\}$. Let $\mathfrak{M}\in \mathcal{D}_\clubsuit$ 
    and $w\in\mathfrak{M}$. In order to prove $\,\vdash_\clubsuit\, \subseteq\, \,\vDash_{\mathfrak{M},w}$, it suffices to show that $\,\vDash_{\mathfrak{M},w}$ satisfies all the rules defining $\,\vdash_\clubsuit$. This readily follows from the previous two lemmas.
\end{proof}


\subsection{Completeness}
To prove completeness with respect to Kripke semantics, the standard approach involves constructing a canonical model with special sets of formulas; each set intuitively gathers precisely the formulas true at this very point. 

Usually, when constructing the canonical model for a non-classical propositional logic with local disjunction (such as \textbf{IPL}, \textbf{BPL} and some subintuitionistic logics discussed in \cite{closer}), the chosen sets are asked to be consistent, deduction closed, and disjunction complete -- for any formula $\varphi,\psi$, if $\varphi\lor\psi\in\Gamma$, then either $\varphi\in\Gamma$ or $\psi\in\Gamma$. 

However, in a restricted language like here, requiring disjunction completeness is inappropriate, due to the absence of local disjunction. Yet we do have our way of expressing the local disjunction -- expression "$\phi \vdash (\alpha_1, \ldots, \alpha_n)$" as previously defined. Thus, the appropriate "disjunction completeness" in our settings should be : whenever $\phi \vdash (\alpha_1, \ldots, \alpha_n)$, if $\phi\in\Gamma$, then $\alpha_1\in\Gamma$ or ... or $\alpha_n\in\Gamma$.

As to the accessibility relation within the canonical model, intuitively it should be the strongest relation $R_c \subseteq \wp(Form_\neg)^2$ expressed with syntactic notions, satisfying that for any $\mathfrak{M}\in \mathcal{D}$, and any $w,v\in\mathfrak{M}$, if $v$ is accessible from $w$ in $\mathfrak{M}$, then $Th(\mathfrak{M},w)\, R_c\, Th(\mathfrak{M},w)$ (- where $\mathcal{D}$ is the class of models under consideration and $Th(\mathfrak{M},w)$ collects all the formulas true at $\mathfrak{M},w$). 

\begin{definition}[\(\,\vdash\)-consistent/closed/prime]
    Let $\ \,\vdash\ \subseteq \wp(Form_\neg) \times Form_\neg$. For any \(\Gamma \subseteq Form_\neg\) and any $\phi\in Form_\neg$, define:
 \begin{itemize}
    \item \(\Gamma\) is \(\,\vdash\)-consistent \(\iff\) \(\Gamma \,\nvdash \bot\);

    \item \(\Gamma\) is \(\,\vdash\)-\(\phi\)-maximal consistent \(\iff\)  \(\Gamma \,\nvdash \phi\) and for any \( \alpha \in \text{Form}_\neg \setminus \Gamma\), \(\Gamma \cup \{\alpha\} \,\vdash\phi\);
    
    \item \(\Gamma\) is \(\,\vdash\)-closed  \(\iff\) for any \( \alpha \in \text{Form}_\neg, (\Gamma \,\vdash \alpha \Rightarrow \alpha \in \Gamma)\);

    \item \(\Gamma\) is \(\,\vdash\)-prime \(\iff\) \ for any \( 2\leq n \in \mathbb{N} \), and for any \( \phi, \alpha_1, \ldots, \alpha_n\in Form_\neg \), if \(\phi \in \Gamma \) and \(\phi \,\vdash (\alpha_1, \ldots, \alpha_n)\), then \(\alpha_1 \in \Gamma \) or \(\ldots \) or \(\alpha_n \in \Gamma\).
 \end{itemize}
\end{definition}

The next lemma shows that being \(\,\vdash\)-\(\phi\)-maximal consistent implies being \(\,\vdash\)-consistent, \(\,\vdash\)-closed and \(\,\vdash\)-prime.

\begin{lemma}\label{max.con. imply con. closed prime}
Let $\ \,\vdash\ \subseteq \wp(Form_\neg) \times Form_\neg$ satisfy \( (\bot)\) and \((Cut) \). For any \( \Gamma \subseteq \text{Form}_\neg \), for any \( \phi \in \text{Form}_\neg \), if \(\Gamma\) is \(\,\vdash\)-\(\phi\)-maximal consistent, then  \(\Gamma\) is  \(\,\vdash\)-consistent, \(\,\vdash\)-closed and \(\,\vdash\)-prime.
\end{lemma}

\begin{proof}
Let \(\Gamma\) be \(\,\vdash\)-\(\phi\)-maximal consistent.

By \( (\bot)\) and \((Cut) \), it is easy to prove that \( \Gamma \nvdash \bot \). 

Let $\alpha$ be an arbitrary formula and \( \Gamma \vdash \alpha \). Suppose
(towards contradiction) that \( \alpha \notin \Gamma \). Then, by the maximality of \( \Gamma \),  \( \Gamma \cup \{ \alpha \} \vdash \phi \). By \( (Cut) \), we have \( \Gamma \vdash \phi \). Contradiction.

Let \( 2 \leq n \in \mathbb{N}^* \), and \( \beta, \alpha_1, \ldots, \alpha_n \) be formulas such that \( \beta \in \Gamma \) and \( \beta \vdash (\alpha_1, \ldots, \alpha_n) \). Suppose
(towards contradiction) that none of \( \alpha_1, \ldots, \alpha_n\) belongs to \(\Gamma \). Then, by the maximality of \( \Gamma \), for any \( i \in \{ 1, \ldots, n \}, \Gamma \cup \{ \alpha_i \} \vdash \phi \). Thus, combining \( \beta \vdash (\alpha_1, \ldots, \alpha_n) \) we get \( \Gamma \cup \{ \beta \} \vdash \phi \). Since \( \beta \in \Gamma \), we have \( \Gamma \vdash \phi \). Contradiction.
\end{proof}

\begin{definition}[\(W^c_{\,\vdash},\ V^c_{\,\vdash}\)]
Let $\ \,\vdash\ \subseteq \wp(Form_\neg) \times Form_\neg$. Define:
\begin{enumerate}[(i)]
    \item \(W^c_{\,\vdash} = \{ \Gamma \subseteq \text{Form}_\neg \mid \Gamma\) is \(\vdash\)-closed,  $\vdash$-consistent and  $\,\vdash$-prime\(\},\)
    \item \(V^c_{ \,\vdash}\) is a function from  \(W^c_{\,\vdash}\) to \(\mathcal{P}(PL)\) such that \(V^c_{ \,\vdash}(\Gamma) = \{ p \in PL \mid p \in \Gamma \}\) for any $\Gamma\in W^c_{\,\vdash}$.
\end{enumerate}
\end{definition}

\begin{definition}[$R^c_{\,\vdash_\clubsuit}$]
Let $\neg[\Gamma]$ denote $\{\neg\alpha\mid\alpha\in\Gamma\}$. Define:
\begin{itemize}
    \item \(R^c_{ \,\vdash_{\mathbf{K}^p}} = \{ (\Gamma, \Delta) \in (W^c_{ \,\vdash_{\mathbf{K}^p}})^2 \mid \Gamma \cap \neg[\Delta] = \emptyset\), and \(\neg T\in\Gamma\) implies \(\Gamma\subseteq\Delta\}\),
    \item \(R^c_{ \,\vdash_{\mathbf{V}}} = \{ (\Gamma, \Delta) \in (W^c_{ \,\vdash_{\mathbf{V}}})^2 \mid \Gamma \cap \neg[\Delta] = \emptyset\), and \(\Gamma\subseteq\Delta\}\),
    \item \(R^c_{ \,\vdash_{\mathbf{B}^p}} = \{ (\Gamma, \Delta) \in (W^c_{ \,\vdash_{\mathbf{B}^p}})^2 \mid \Gamma \cap \neg[\Delta] = \emptyset\}\),
    \item \(R^c_{ \,\vdash_{\mathbf{B4}^p}} = \{ (\Gamma, \Delta) \in (W^c_{ \,\vdash_{\mathbf{B4}^p}})^2 \mid \Gamma \cap \neg[\Delta] = \emptyset\}\).
\end{itemize}
\end{definition}

\begin{remark}
    Let $R'^c_{ \,\vdash_{\mathbf{B}^p}} = \{ (\Gamma, \Delta) \in (W^c_{ \,\vdash_{\mathbf{B}^p}})^2 \mid \Gamma \cap \neg[\Delta] = \emptyset$, and for any $\alpha\in Form_\neg$, $(\neg\neg\alpha\in\Gamma\cap\Delta\Rightarrow\alpha\in\Delta)\}$, and let \(R'^c_{ \,\vdash_{\mathbf{B4}^p}} = \{ (\Gamma, \Delta) \in (W^c_{ \,\vdash_{\mathbf{B4}^p}})^2 \mid \Gamma \cap \neg[\Delta] = \emptyset\), and \(\Gamma\subseteq\Delta\}\). 
    In fact, we can demonstrate that $R^c_{ \,\vdash_{\mathbf{B}^p}} = R'^c_{ \,\vdash_{\mathbf{B}^p}}$ and $R^c_{ \,\vdash_{\mathbf{B4}^p}} = R'^c_{ \,\vdash_{\mathbf{B4}^p}}$, using the primeness of elements in $W^c_{ \,\vdash_{\clubsuit}}$, together with the fact that $\,\vdash_{\mathbf{B}^p}$ satisfies (Prop${}_{\mathbf{B}^p}$)+($\neg\neg$I) and $\,\vdash_{\mathbf{B4}^p}$ satisfies (Prop${}_{\mathbf{B4}^p}$).  
\end{remark}

\begin{definition}[canonical model for $\vdash_{\clubsuit}$]
For any $\clubsuit\in\{\mathbf{K}^p, \mathbf{V}, \mathbf{B}^p, \mathbf{B4}^p\}$, define:
\[
\mathcal{M}^c_{ \vdash_\clubsuit} = ( \mathcal{F}^c_{ \,\vdash_\clubsuit}, V^c_{ \,\vdash_\clubsuit} ) = ( W^c_{ \,\vdash_\clubsuit}, R^c_{ \,\vdash_\clubsuit}, V^c_{ \,\vdash_\clubsuit} ).
\]
\end{definition}

Next, we present several crucial lemmas for the completeness theorem.

\begin{lemma}[Lindenbaum's lemma]
Let $\ \,\vdash\ \subseteq \wp(Form_\neg) \times Form_\neg$ satisfy \( (A), (Cut)\), \((Com), (Mon) \).
If \(\Gamma \nvdash\phi\), then there exists \(\Phi \subseteq \text{Form}_\neg\), such that \(\Phi\) is \(\,\vdash\)-\(\phi\)-maximal consistent and \(\Gamma \subseteq \Phi\) and $\phi\notin\Phi$.
\end{lemma}

\begin{proof}
Consider $\Omega = \{ \Delta \subseteq \text{Form}_\neg \mid \Gamma \subseteq \Delta, \Delta \text{ is } \vdash\text{-closed, } \phi\notin\Delta \}$. 

It is easy to show that $\text{Th}_{\,\vdash}(\Gamma) =\{\alpha\mid\Gamma\vdash\alpha\}$ is in $\Omega$ (using (A), (Com), (Cut), and $\Gamma\nvdash\phi$), and every chain in $\Omega$ has an upper bound (using (Com), (Mon)). 

By Zorn's lemma, $\Omega$ has a maximal element, denoted as $\Phi$. We now show that $\Phi$ is $\vdash$-$\phi$-maximal consistent. 

Since $\phi\notin\Delta$ and $\Phi$ is $\vdash$-closed, we have $\Phi\nvdash\phi$. 

For any $\psi \in \text{Form}_\neg \setminus \Phi$, since $\text{Th}_{\,\vdash}(\Phi \cup \{\psi\}) \supset \Phi \supseteq \Gamma$ and it is $\vdash$-closed (using (Com), (Cut)), by the maximality of $\Phi$, we have $\phi\in\text{Th}_{\,\vdash}(\Phi \cup \{\psi\})$. So $\Phi \cup \{\psi\}\vdash\phi$.
\end{proof}

\begin{lemma}[Existence lemma]
For any $\clubsuit\in\{\mathbf{K}^p, \mathbf{V}, \mathbf{B}^p, \mathbf{B4}^p\}$, for any $\Gamma \in W^c_{ \,\vdash_\clubsuit}$ and $\alpha \in Form_\neg$, if $\neg \alpha \notin \Gamma$, then there exists $ \Phi \in W^c_{ \,\vdash_\clubsuit}$, such that $\Gamma R^c_{ \,\vdash_\clubsuit} \Phi$  and $\alpha\in \Phi. $
\end{lemma}

\begin{proof}
Consider $\Omega = \{ \Delta \mid \Delta$ is $\vdash_\clubsuit$-consistent and $\vdash_\clubsuit$-closed, $\Gamma \cap \neg[\Delta] = \emptyset$ and $\alpha\in \Delta \}$.

Note that $Th_{\,\vdash_\clubsuit}(\{\alpha\}) \in \Omega$ (using ($\neg\bot$), ($\neg$ antitone), (Cut), (Mon), (Com), (A), with $\Gamma$ $\,\vdash_\clubsuit$-closed and $\neg \alpha \notin \Gamma$). So $\Omega$ is not empty.

It is easy to show (using (Com), (Mon)) that for any directed family $\mathcal{H} \subseteq \Omega$, $\bigcup \mathcal{H} \in \Omega$.

By Zorn's lemma, $\Omega$ has a maximal element, denoted as $\Phi$. We prove that $\Phi$ meets all the requirements.

\begin{enumerate}[(i)]
    \item\label{i} \textbf{Task 1:} $\Phi \in W^c_{ \,\vdash_\clubsuit}$.
    
It suffices to show $\Phi$ is $\,\vdash_\clubsuit$-prime.

Let $2 \leq n \in \mathbb{N}^*$, and $\beta, \alpha_1, \ldots, \alpha_n\in Form_\neg$. Assume that $\beta \in \Phi$ and $\beta \vdash_{\clubsuit} (\alpha_1, \ldots, \alpha_n)$.

Suppose towards contradiction that none of $\alpha_1, \ldots, \alpha_n$ belongs to $\Phi$.

Consider $\text{Th}_{\,\vdash_\clubsuit}(\Phi \cup \{\alpha_i\})$.

Note that $\Phi\subset\text{Th}_{\,\vdash_\clubsuit}(\Phi \cup \{\alpha_i\})$ and $\text{Th}_{\,\vdash_\clubsuit}(\Phi \cup \{\alpha_i\})$ is $\,\vdash_\clubsuit$-closed ((Com), (Cut), (Mon)), and ($\bot \in \Delta \Rightarrow \Gamma \cap \neg[\Delta] \neq \emptyset$) ($\Gamma$  $\,\vdash_\clubsuit$-closed, ($\neg\bot$), (Mon)). By the maximality of $\Phi$, we have $\Gamma \cap \neg[\text{Th}_{\,\vdash_\clubsuit}(\Phi \cup \{\alpha_i\})] \neq \emptyset$. Namely, for any $i \in \{1, \ldots, n\}$, there exists $\neg\chi_i\in\Gamma$ with $\Phi \cup \{\alpha_i\}\,\vdash_\clubsuit\chi_i$.

Using (Com), (Mon), ($\land$I/E), (Cut) and closedness of $\Phi$, we can infer that there exists $\psi \in \Phi$, for any $i \in \{1, \ldots, n\}$,  $\psi \land \alpha_i\,\vdash_\clubsuit\chi_i$.

Using ($\neg$ antitone), we have $\neg\chi_i\,\vdash_\clubsuit\neg(\psi \land \alpha_i)$ for any $i \in \{1, \ldots, n\}$. Since $\neg\chi_i\in\Gamma$ and $\Gamma$ is $\,\vdash_\clubsuit$-closed, we get $\neg(\psi \land \alpha_i) \in \Gamma$ for any $i \in \{1, \ldots, n\}$.

Moreover, since $\beta \vdash_{\clubsuit} (\alpha_1, \ldots, \alpha_n)$, we have $\neg(\psi \land \alpha_1) \land \ldots \land \neg(\psi \land \alpha_n) \vdash_{\clubsuit} \neg(\psi \land \beta)$.

Thus, using (A), ($\land$I), (Cut) and closedness of $\Gamma$, we have $\neg(\psi \land \beta) \in \Gamma$.

On the other hand, since $\psi, \beta \in \Phi$, using ($\land$I) and closedness of $\Phi$, we have $\psi \land \beta \in \Phi$, contradicting $\Gamma \cap \neg[\Phi] = \emptyset$.

  \item \textbf{Task 2:} $\Gamma R^c_{ \,\vdash_\clubsuit} \Phi$.

It breaks down into three cases.

    \begin{enumerate}[(a)]
        \item $\clubsuit=\mathbf{K}^p$.

              Then by definition, $\Gamma R^c_{ \,\vdash_\clubsuit} \Phi\iff$\(\Gamma \cap \neg[\Phi] = \emptyset\) and \((\neg T\in\Gamma\Rightarrow\Gamma\subseteq\Phi)\).
              
              Firstly, $\Gamma \cap \downarrow\neg[\Phi] = \emptyset$ is true by $\Phi\in\Omega$.

              Secondly, assuming $\neg T \in \Phi$, we prove $\Gamma\subseteq\Phi$.
              
              Let $\phi \in \Gamma$ be arbitrary and suppose (towards contradiction) that $\phi \notin\Phi$.
              
              Using similar reasoning as in (\ref{i}), it can be shown that there exists $\psi \in \Phi$ with $\neg(\psi \land \phi) \in \Gamma$.
              
              By (Prop${}_{\mathbf{K}^p}$) (since $\clubsuit=\mathbf{K}^p$), $\phi \land \neg(\psi \land \phi) \vdash_\clubsuit \neg(\psi \land \neg T)$.

              Since $\phi \in \Gamma$ and $\neg(\psi \land \phi) \in \Gamma$, by ($\land$I), (Cut) and closedness of $\Gamma$, we have $\neg(\psi \land \neg T) \in \Gamma$.
              
              However, since $\psi, \neg T \in \Phi$, by ($\land$I), (Cut) and closedness of $\Phi$, we have $\psi \land \neg T \in \Phi$, contradicting $\Gamma \cap \neg[\Phi] = \emptyset$.

        \item $\clubsuit=\mathbf{V}$.

              Then by definition, $\Gamma R^c_{ \,\vdash_\clubsuit} \Phi\iff$\(\Gamma \cap \neg[\Phi] = \emptyset\) and \(\Gamma\subseteq\Phi\).
              
              Firstly, $\Gamma \cap \downarrow\neg[\Phi] = \emptyset$ is true by $\Phi\in\Omega$.

              Secondly, we prove $\Gamma\subseteq\Phi$.
              
              Let $\phi \in \Gamma$ be arbitrary and suppose (towards contradiction) that $\phi \notin\Phi$.
              
              Using similar reasoning as in (\ref{i}), it can be shown that there exists $\psi \in \Phi$ with $\neg(\psi \land \phi) \in \Gamma$.
              
              By (Prop${}_{\mathbf{V}}$) (since $\clubsuit=\mathbf{V}$), $\phi \land \neg(\psi \land \phi) \vdash_\clubsuit \neg\psi$.

              Since $\phi \in \Gamma$ and $\neg(\psi \land \phi) \in \Gamma$, by ($\land$I), (Cut) and closedness of $\Gamma$, we have $\neg\psi \in \Gamma$.
              
              However,  $\psi \in \Phi$, contradicting $\Gamma \cap \neg[\Phi] = \emptyset$.

        \item $\clubsuit=\mathbf{B}^p$ or $\clubsuit=\mathbf{B4}^p$.

              Then by definition, $\Gamma R^c_{ \,\vdash_\clubsuit} \Phi\iff$\(\Gamma \cap \neg[\Phi] = \emptyset\), which is true by $\Phi\in\Omega$.        
    \end{enumerate}
\end{enumerate}
\end{proof}

\begin{lemma}[Truth Lemma]
For any $\clubsuit\in\{\mathbf{K}^p, \mathbf{V}, \mathbf{B}^p, \mathbf{B4}^p\}$, for any $\Gamma \in W^c_{ \,\vdash_\clubsuit}$ and $\alpha \in Form_\neg$, 
\[\mathfrak{M}^c_{ \,\vdash_\clubsuit},\Gamma\,\vDash\alpha \iff \alpha\in\Gamma.\]
\end{lemma}

\begin{proof}
Let $\clubsuit\in\{\mathbf{K}^p, \mathbf{V}, \mathbf{B}^p, \mathbf{B4}^p\}$. Use induction on formulas to show that for any $\Gamma \in W^c_{ \,\vdash_\clubsuit}$, the target equivalence holds. The case for propositional letters follows from the definition of $V^c_{\,\vdash_\clubsuit}$. For any $\Gamma\in W^c_{\,\vdash_\clubsuit}$, $\bot\notin\Gamma$, since $\Gamma$ is ${\,\vdash_\clubsuit}$-consistent and $\,\vdash_\clubsuit$ satisfies (A). The case for $\land$ is proved using ($\land$I), ($\land$E) and the closedness of each $\Gamma\in W^c_{\,\vdash_\clubsuit}$. The case for $\neg$ follows from the previous lemma and the fact that \(\Gamma \cap \neg[\Delta] = \emptyset\) for any $(\Gamma,\Delta)\in R^c_{\,\vdash_\clubsuit}$.
\end{proof}

\begin{lemma}\label{canoical model belongs to D}
   For any $\clubsuit\in\{\mathbf{V}, \mathbf{B}^p, \mathbf{B4}^p\}$, $\mathfrak{M}^c_{ \,\vdash_\clubsuit}\in \mathcal{D}_\clubsuit$. Moreover, $\mathfrak{M}^c_{ \,\vdash_{\mathbf{K}^p}}\in \mathcal{D}_{\mathbf{K}^p-}$.
\end{lemma}
\begin{proof}
In this proof, we omit the subscripts from $\Vert\cdot\Vert_{\mathfrak{M}^c_{ \,\vdash_{\clubsuit}}}$ and from set operators such as $\Box_{\mathfrak{F}^c_{ \,\vdash_{\clubsuit}}}$. 
\begin{enumerate}[(1)]
    \item 
    \indent\quad $\mathfrak{M}^c_{ \,\vdash_{\mathbf{K}^p}}\in \mathcal{D}_{\mathbf{K}^p-}$ 
    
    $\Leftrightarrow$ for any $\alpha\in Form_{\neg}$,  $R^c_{ \,\vdash_{\mathbf{K}^p}}\Vert\alpha\Vert \cap\Vert\neg T\Vert \subseteq \Vert\alpha\Vert$  \hfill (definition of $\mathcal{D}_{\mathbf{K}^p-}$, $\Vert\neg T\Vert=\Box\emptyset$)
        
    $\Leftrightarrow$ for any $\alpha\in Form_{\neg}$, for any $\Gamma,\Delta\in W^c_{ \,\vdash_{\mathbf{K}^p}}$, if $\alpha\in\Gamma$ and $\Gamma R^c_{ \,\vdash_{\mathbf{K}^p}} \Delta$ and $\neg T\in\Delta$, then $\alpha\in\Delta$ 
    \hfill \quad(definition of $\,\vDash$, Truth lemma)
    
    --which is true by the definition of $R^c_{ \,\vdash_{\mathbf{K}^p}}$.

    \item  First, we prove that $R^c_{ \,\vdash_{\mathbf{V}}}$ is transitive. Assume $\Gamma R^c_{ \,\vdash_{\mathbf{V}}}\Delta$ and $\Delta R^c_{ \,\vdash_{\mathbf{V}}} \Phi$. By definition, \(\Gamma \cap \neg[\Delta] = \emptyset\) and \(\Gamma\subseteq\Delta\), and \(\Delta \cap \neg[\Phi] = \emptyset\) and \(\Delta\subseteq\Phi\). Then on the one hand, from \(\Gamma\subseteq\Delta\) and \(\Delta\subseteq\Phi\) we have \(\Gamma\subseteq\Phi\). On the other hand, from \(\Gamma\subseteq\Delta\) and \(\Delta \cap \neg[\Phi] = \emptyset\) we have \(\Gamma\cap \neg[\Phi] = \emptyset\). So, $\Gamma R^c_{ \,\vdash_{\mathbf{V}}}\Phi$.
    
    Then, 
    
    \indent\quad $\mathfrak{M}^c_{ \,\vdash_{\mathbf{V}}}\in \mathcal{D}_{\mathbf{V}}$ 
    
    $\Leftrightarrow$ for any $\alpha\in Form_{\neg}$, $\Vert \alpha \Vert \subseteq \Box\Vert\alpha \Vert$ \hfill (definition of $\mathcal{D}_{\mathbf{V}}$, $R^c_{ \,\vdash_{\mathbf{V}}}$ is transitive, lemma \ref{formula.fixedpoint}, lemma \ref{frame condition.fixpoint})
    
    $\Leftrightarrow$ for any $\alpha\in Form_{\neg}$, for any $\Gamma,\Delta\in W^c_{ \,\vdash_{\mathbf{V}}}$, if $\alpha\in\Gamma$ and $\Gamma R^c_{ \,\vdash_{\mathbf{V}}} \Delta$, then $\alpha\in\Delta$ 
    \hfill \quad(definition of $\,\vDash$, Truth lemma)
    
    --which is true by the definition of $R^c_{ \,\vdash_{\mathbf{V}}}$.

    \item First, we prove that $R^c_{ \,\vdash_{\mathbf{B}^p}}$ is symmetric. Assume $\Gamma R^c_{ \,\vdash_{\mathbf{B}^p}}\Delta$. By definition, \(\Gamma \cap \neg[\Delta] = \emptyset\). Suppose (towards contradiction) that \(\Delta \cap \neg[\Gamma] \neq \emptyset\). Then there is $\alpha\in\Gamma$ such that $\neg\alpha\in\Delta$. Since $\,\vdash_{\mathbf{B}^p}$ satisfies ($\neg\neg$I) and $\Gamma$ is $\,\vdash_{\mathbf{B}^p}$-closed, we have $\neg\neg\alpha\in\Gamma$, which contradicts \(\Gamma \cap \neg[\Delta] = \emptyset\). So, $\Delta R^c_{ \,\vdash_{\mathbf{V}}}\Gamma$.

    Then, 
    
    \indent\quad 
    $\mathfrak{M}^c_{ \,\vdash_{\mathbf{B}^p}}\in \mathcal{D}_{\mathbf{B}^p}$ 
    
    $\Leftrightarrow$ for any $\alpha\in Form_{\neg}$, $\Vert\neg\neg\alpha\Vert \subseteq \Vert \neg\alpha\Vert\cup \Vert\alpha \Vert$ \hfill (definition of $\mathcal{D}_{\mathbf{B}^p}$, $R^c_{ \,\vdash_{\mathbf{B}^p}}$ is symmetric, lemma \ref{formula.fixedpoint}, lemma \ref{frame condition.fixpoint},  $\neg\Vert\alpha\Vert=\Vert\neg\alpha\Vert$)
        
    $\Leftrightarrow$ for any $\alpha\in Form_{\neg}$, for any $\Gamma\in W^c_{ \,\vdash_{\mathbf{B}^p}}$, if $\neg\neg\alpha\in\Gamma$, then $\alpha\in\Gamma$ or  $\neg\alpha\in\Gamma$
    \hfill \quad(definition of $\,\vDash$, Truth lemma)
    
    --which is true since $\,\vdash_{\mathbf{B}^p}$ satisfies (Prop${}_{\mathbf{B}^p}$) and $\Gamma$ is $\,\vdash_{\mathbf{B}^p}$-prime.

    \item First, we prove that $R^c_{ \,\vdash_{\mathbf{B4}^p}}$ is both symmetric and transitive. 
    Symmetry is proved analogously as in (3). As for the transitivity, assume  that $\Gamma R^c_{ \,\vdash_{\mathbf{V}}}\Delta$ and $\Delta R^c_{ \,\vdash_{\mathbf{V}}} \Phi$. By definition, \(\Gamma \cap \neg[\Delta] = \emptyset\) and \(\Delta \cap \neg[\Phi] = \emptyset\). Suppose (towards contradiction) that \(\Gamma \cap \neg[\Phi] \neq \emptyset\). Then there is $\alpha\in\Phi$ such that $\neg\alpha\in\Gamma$. Then $\neg\alpha\notin\Delta$ since $\alpha\in\Phi$ and \(\Delta \cap \neg[\Phi] = \emptyset\). Then, since $\,\vdash_{\mathbf{B4}^p}$ satisfies (Prop${}_{\mathbf{B4}^p}$) and $\Delta$ is $\,\vdash_{\mathbf{B}^p}$-prime, we have $\alpha\in\Delta$. But we already have $\neg\alpha\in\Gamma$, which contradicts \(\Gamma \cap \neg[\Delta] = \emptyset\). So, $\Gamma R^c_{ \,\vdash_{\mathbf{V}}}\Phi$.    

    Then, 
    
    \indent\quad 
    $\mathfrak{M}^c_{ \,\vdash_{\mathbf{B4}^p}}\in \mathcal{D}_{\mathbf{B4}^p}$ 
    
    $\Leftrightarrow$ for any $\alpha\in Form_{\neg}$, $W^c_{ \,\vdash_{\mathbf{B4}^p}}= \Vert \neg\alpha\Vert\cup \Vert\alpha \Vert$ \hfill (definition of $\mathcal{D}_{\mathbf{B4}^p}$, $R^c_{ \,\vdash_{\mathbf{B}^p}}$ is symmetric and transitive, lemma \ref{formula.fixedpoint}, lemma \ref{frame condition.fixpoint})
        
    $\Leftrightarrow$ for any $\alpha\in Form_{\neg}$, for any $\Gamma\in W^c_{ \,\vdash_{\mathbf{B4}^p}}$,  $\alpha\in\Gamma$ or  $\neg\alpha\in\Gamma$
    \hfill \quad(definition of $\,\vDash$, Truth lemma)
    
    --which is true since $\,\vdash_{\mathbf{B4}^p}$ satisfies (Prop${}_{\mathbf{B4}^p}$) and $\Gamma$ is $\,\vdash_{\mathbf{B4}^p}$-prime.
   
\end{enumerate}
\end{proof}

\begin{theorem}[Completeness]\label{completeness part 1}
For each $\clubsuit\in\{\mathbf{K}^p, \mathbf{V}, \mathbf{B}^p, \mathbf{B4}^p\}$, for any $\Gamma\subseteq Form_\neg$ and $\phi\in Form_\neg$,
    \[\Gamma\,\vDash_\clubsuit\phi \Rightarrow\Gamma\,\vdash_\clubsuit\phi.\]
\end{theorem}
\begin{proof}
    Let $\clubsuit\in\{\mathbf{K}^p, \mathbf{V}, \mathbf{B}^p, \mathbf{B4}^p\}$. Let $\Gamma\subseteq Form_\neg$ and $\phi\in Form_\neg$. We prove the contrapositive of the target entailment. 
    
    Assume that $\Gamma\,\nvdash_\clubsuit\phi$. Then by Lindenbaum's lemma and lemma \ref{max.con. imply con. closed prime}, there exists \(\Phi \in W^c_{ \,\vdash_\clubsuit}\) with \(\Gamma \subseteq \Phi\) and $\phi\notin\Phi$. By the Truth lemma, $\mathfrak{M}^c_{ \,\vdash_\clubsuit},\Phi\,\vDash\Gamma$ and $\mathfrak{M}^c_{ \,\vdash_\clubsuit},\Phi\,\nvDash\alpha$. Moreover, according to lemma \ref{canoical model belongs to D}, $\mathfrak{M}^c_{ \,\vdash_\clubsuit}\in \mathcal{D}_\clubsuit$ for $\clubsuit\in\{\mathbf{V}, \mathbf{B}^p, \mathbf{B4}^p\}$ while $\mathfrak{M}^c_{ \,\vdash_{\mathbf{K}^p}}\in \mathcal{D}_{\mathbf{K}^p-}$. So $\Gamma\,\nvDash_\clubsuit\phi$ if $\clubsuit\in\{\mathbf{V}, \mathbf{B}^p, \mathbf{B4}^p\}$, and $\Gamma\,\nvDash_{\mathbf{K}^p-}\phi$ if $\clubsuit=\mathbf{K}^p$. Finally, by lemma \ref{vDash_K-=vDash_K, vDash_T-=vDash_T}, we have $\Gamma\,\nvDash_{\mathbf{K}^p}\phi$ if $\clubsuit=\mathbf{K}^p$.
\end{proof}


\section{Axiomatization in $\{\bot,\land,\to\}$-language}

In this section, we only consider the formulas built up from propositional letters and $\bot$ using $\{\land,\to\}$.
Let $Form_\to$ collect all these formulas. Throughout this section, let us make the following conventions. 
\begin{enumerate}[(1)]
    \item $\neg \alpha := \alpha\to\bot$;
    \item $T := \neg\bot$;
    \item $\Box \alpha := T\to\bot$. 
\end{enumerate}
    
The concept of ``$\phi \vdash (\alpha_1, \ldots, \alpha_n)$" needs to be updated due to the increase in expressive power of the language. 

\begin{definition}[``$\phi \vdash (\alpha_1, \ldots, \alpha_n)$"]
    For any $\ \vdash\ \subseteq \wp(Form_\to) \times Form_\to$, for any \( 2\leq n \in \mathbb{N} \), and for any \( \phi, \alpha_1, \ldots, \alpha_n\in Form_\to \), we denote ``$\phi \vdash (\alpha_1, \ldots, \alpha_n)$" if and only if :
    \begin{enumerate}[(i)]
        \item for any $\Gamma\subseteq Form_\to$ and $\psi\in Form_\to$, if $\Gamma \cup \{ \alpha_i \} \vdash \psi$ for any $i \in \{ 1, \ldots, n \}$, then $\Gamma \cup \{ \phi \} \vdash \psi$;
        \item for any $\psi,\beta\in Form_\to$, $(\psi \land \alpha_1\to\beta) \land \ldots \land (\psi \land \alpha_n\to\beta) \vdash (\psi \land \phi\to\beta)$.
    \end{enumerate}
\end{definition}

Again, we can define $\,\vDash_{\mathfrak{M},w} \subseteq \wp(Form_\to) \times Form_\to$ and check that:
\begin{lemma}\label{pseudo disjunction.soundness'}
    For any Kripke model $\mathfrak{M}$, for any \( 2\leq n \in \mathbb{N} \), and for any \( \phi, \alpha_1\),...,\(\alpha_n\in Form_\to \), if $\Vert\phi\Vert_{\mathfrak{M}} \subseteq \Vert \alpha_1\Vert_{\mathfrak{M}} \cup ... \cup\Vert\alpha_n \Vert_{\mathfrak{M}}$, then for any $w\in \mathfrak{M}$, $\phi \,\vDash_{\mathfrak{M},w} (\alpha_1, \ldots, \alpha_n)$. 
\end{lemma}

\subsection{Sequent calculi}
Recall the rules (A),(Cut),(\(\land\)I),(\(\land\)E) and (\(\bot\)) defined in Section 3.1. In addition, we introduce some new ones:

\begin{itemize}    
    \item[] (DT${}_0$)\quad $\alpha \vdash \beta\  \Rightarrow \ \ \vdash\alpha\to\beta$ ;
    \item[] ($\to\land$)\quad $(\alpha\to\beta)\land(\alpha\to\gamma)\,\vdash\alpha\to \beta\land\gamma$ ;
    \item[] ($\to$tr)\quad $(\alpha\to\beta)\land(\beta\to\gamma)\,\vdash\alpha\to \gamma$ ;
    
    \item[](Refl)\quad $\alpha\land(\alpha\to\beta)\vdash\beta$ ;
    \item[](Sym)\quad $\alpha\vdash (\beta,\neg(\alpha\to\beta))$ ;

    \item[](Prop${}_{\mathbf{K}^p-}$)\quad $\alpha\vdash \neg T\to \alpha$ ;
    \item[](Prop${}_{tr}$)\quad $\alpha\vdash \Box \alpha$ ;
    \item[](Prop${}_{sy}$)\quad $\alpha\vdash \neg\neg \alpha \to \alpha$ .
\end{itemize}

    

Now, the rules (A),(Cut),(\(\land\)I),(\(\land\)E),(\(\bot\)),(DT${}_0$),($\to\land$) and ($\to$tr) are called basic rules. 

\begin{definition}
  Let \( \vdash_{\mathbf{K}^p} \) (and $\vdash_{\mathbf{T}^p}$, $\vdash_{\mathbf{V}}$, $\vdash_{\mathbf{B}^p}$, $\vdash_{\mathbf{B4}^p}$, $\vdash_{\mathbf{I}}$, $\vdash_{\mathbf{O}}$, $\vdash_{\mathbf{C}}$, respectively) be the smallest $\ \vdash\ \subseteq \wp(Form_\to) \times Form_\to$ that satisfies, in addition to the basic rules, (Prop${}_{\mathbf{K}^p-}$) (and (Refl), (Prop${}_{tr}$), 
(Sym)+(Prop${}_{sy}$), (Sym)+(Prop${}_{tr}$), (Refl)+(Prop${}_{tr}$), (Refl)+(Prop${}_{sy}$), (Refl)+(Sym)+(Prop${}_{tr}$), respectively).
\end{definition}

It is easy to check that:
\begin{lemma}
    for each $\clubsuit\in\{\mathbf{K}^p, \mathbf{T}^p, \mathbf{B}^p, \mathbf{V}, \mathbf{B4}^p, \mathbf{I}, \mathbf{O}, \mathbf{C}\}$, $\vdash_{\clubsuit}$ enjoys (Com).
\end{lemma}

\subsection{Soundness}
Recall that a rule $(X)$ is correct at a pointed model $\mathfrak{M},w$ if and only if $\ \vDash_{\mathfrak{M},w}$ satisfies $(X)$. 
It is easy to show that :
\begin{lemma}
    All basic rules, except (DT${}_0$), are correct at any pointed Kripke model.
\end{lemma}

For (DT${}_0$), we need the notion of correctness at a model (instead of a pointed model).

\begin{definition}[$\,\vDash_{\mathfrak{M}}$]
    Let $\mathfrak{M}$ be an arbitrary Kripke model. Define $\,\vDash_{\mathfrak{M}}\ \subseteq\ $ $\wp(Form_\to)\times Form_\to$ to be: for any $\Gamma\subseteq Form_\to$ and $\phi\in Form_\to$,
    \begin{align*}
       \Gamma\vDash_{\mathfrak{M}}\phi\; &\iff\; \textit{ for any }w\in\mathfrak{M},\ \Gamma \vDash_{\mathfrak{M},w}\phi\\
       &\iff\; \textit{ for any }w\in\mathfrak{M},\ \mathfrak{M},w\vDash\Gamma \textit{ implies } \mathfrak{M},w\vDash\phi\ .
    \end{align*}
\end{definition}

\begin{definition}[Correctness at a model]
    A sequent rule $(X)$ is correct at a Kripke model $\mathfrak{M}$, if and only if $\ \vDash_{\mathfrak{M}}$ satisfies $(X)$.   
\end{definition}

It is not hard to show the following relation between two kinds of correctness. 
\begin{lemma}
    For any Kripke model $\,\mathfrak{M}$, if a rule $(X)$ is correct at $\,\mathfrak{M},w$ for each $w \in\mathfrak{M}$, then $(X)$ is correct at $\,\mathfrak{M}$. The reverse is also true, given that $(X)$ is premise-free, that is, $(X)$ is a axiom.
\end{lemma}

Now we prove the correctness of (DT${}_0$).
\begin{lemma}
    (DT${}_0$) is correct at any Kripke model.
\end{lemma}
\begin{proof}
    Let $\,\mathfrak{M}$ be a Kripke model. We want to show that $\alpha \vDash_{\mathfrak{M}} \beta\  \Rightarrow \ \ \vDash_{\mathfrak{M}}\alpha\to\beta$, for any formula $\alpha$ and $\beta$. Assume $\alpha \vDash_{\mathfrak{M}} \beta$. Then we readily know that for any $w\in\mathfrak{M}$, $\,\vDash_{\mathfrak{M},w}\alpha\to\beta$.
\end{proof}

Next, we examine the correctness of non-basic rules.
\begin{lemma}\label{correctness of re,sy,tr rules}
~
\begin{enumerate}[(1)]
    \item (Refl) is correct at any pointed reflexive model.
    \item (Sym) is correct at any pointed symmetric model.
\end{enumerate}
\end{lemma}

\begin{proof}
~
\begin{enumerate}[(1)]   
    \item Let $\mathfrak{M},w$ be a pointed reflexive model, and let $\alpha,\beta$ be formulas. We want to show that $\alpha\land(\alpha\to\beta)\,\vDash_{\mathfrak{M},w} \beta$. Assume $\mathfrak{M},w\vDash \alpha\land(\alpha\to\beta)$. Then by the definition of satisfaction relation and the reflexivity of $\mathfrak{M}$, we have $\mathfrak{M},w\vDash \alpha$ and ($\mathfrak{M},w\vDash \alpha\Rightarrow \mathfrak{M},w\vDash \beta$). So $\mathfrak{M},w\vDash \beta$. 
    
    \item Let $\mathfrak{M}=(\mathfrak{F},V)=(W,R,V)$ be a symmetric model, and let $\alpha,\beta$ be formulas. We want to show that $\alpha\,\vDash_{\mathfrak{M},w} (\beta,\neg(\alpha\to\beta))$ for any $w\in\mathfrak{M}$. By lemma \ref{pseudo disjunction.soundness'}, it suffices to prove $\Vert\alpha\Vert_{\mathfrak{M}} \subseteq \Vert \beta\Vert_{\mathfrak{M}}\cup\Vert\neg(\alpha\to\beta)\Vert_{\mathfrak{M}}$. By the definition of truth sets and symmetry of $R$, it is equivalent to $\Vert\alpha\Vert_{\mathfrak{M}}\setminus\Vert \beta\Vert_{\mathfrak{M}} \subseteq   \Box_\mathfrak{F}R(\Vert\alpha\Vert_{\mathfrak{M}}\setminus\Vert\beta\Vert_{\mathfrak{M}})$, which is obviusly true.
     
\end{enumerate}
\end{proof}

\begin{lemma}\label{definability}
~
\begin{enumerate}[(1)]
    \item For any model $\mathfrak{M}$, $\ \mathfrak{M}\in \mathcal{D}_{\mathbf{K}^p-}$ 
    $\Leftrightarrow$ (Prop${}_{\mathbf{K}^p-}$) is correct at $\mathfrak{M}$.
    \item For any symmetric model $\mathfrak{M}$, $\ \mathfrak{M}$ is a BIO-model 
    $\Leftrightarrow$ (Prop${}_{sy}$) is correct at $\mathfrak{M}$.  
    \item For any transitive model $\mathfrak{M}$, $\ \mathfrak{M}$ is a BIO-model 
    $\Leftrightarrow$ (Prop${}_{tr}$) is correct at $\mathfrak{M}$.  
\end{enumerate}

\end{lemma}

\begin{proof}
In this proof, we omit the subscript "$\mathfrak{M}$" of $\Vert\cdot\Vert$. 
\begin{enumerate}[(1)]
    \item For any model $\mathfrak{M}=(\mathfrak{F},V)=(W,R,V)$,\\
    \indent\quad $\mathfrak{M}\in \mathcal{D}_{\mathbf{K}^p-}$ 
    
    $\Leftrightarrow$ for any $\alpha\in Form_{\to}$,  $R\Vert\alpha\Vert \cap\Box_{\mathfrak{F}}\emptyset \subseteq \Vert\alpha\Vert$  \hfill (definition of $\mathcal{D}_{\mathbf{K}^p-}$, lemma \ref{formula.fixedpoint.K-})
    
    $\Leftrightarrow$ for any $\alpha\in Form_{\to}$, $\Vert \alpha \Vert \subseteq \Vert \neg T\Vert\to_{\mathfrak{F}} \Vert\alpha \Vert$ \hfill ($RX\subseteq Y\Leftrightarrow X\subseteq \Box_{\mathfrak{F}} Y$, $\Vert\neg T\Vert=\Box_{\mathfrak{F}}\emptyset$)
    
    $\Leftrightarrow$ for any $\alpha\in Form_{\to}$, $\Vert \alpha \Vert \subseteq \Vert \neg T\to\alpha \Vert$  \hfill \quad(definition of $\Vert\cdot\Vert$)
    
    $\Leftrightarrow$ for any $\alpha\in Form_{\to}$, $ \alpha  \,\vDash_{\mathfrak{M}}  \neg T\to\alpha$  \hfill \quad(definition of the satisfaction relation and $\,\vDash_{\mathfrak{M}}$)
    
    $\Leftrightarrow$ (Prop${}_{\mathbf{K}^p-}$) correct at $\mathfrak{M}$  \hfill \quad(definition of correctness).

    \item For any symmetric model $\mathfrak{M}=(\mathfrak{F},V)=(W,R,V)$,\\
    \indent\quad $\mathfrak{M}$ is a BIO-model
    
    $\Leftrightarrow$ for any $\alpha\in Form_{\to}$, $\Vert\alpha \Vert \subseteq \neg_{\mathfrak{F}}\neg_{\mathfrak{F}}\Vert\alpha\Vert \to_{\mathfrak{F}} \Vert\alpha \Vert$ \hfill (definition of BIO-model, lemma \ref{formula.fixedpoint}, $R$ is symmetric, lemma \ref{frame condition.fixpoint})
    
    $\Leftrightarrow$ for any $\alpha\in Form_{\to}$, $\Vert\alpha \Vert \subseteq \Vert\neg\neg\alpha\to\alpha \Vert$  \hfill \quad(definition of $\Vert\cdot\Vert$)
    
    $\Leftrightarrow$ for any $\alpha\in Form_{\to}$, $\neg\neg\alpha \,\vDash_{\mathfrak{M}} (\alpha, \neg\alpha)$  \hfill \quad(definition of the satisfaction relation and $\,\vDash_{\mathfrak{M}}$)
    
    $\Leftrightarrow$ (Prop${}_{sy}$) is correct at $\mathfrak{M}$  \hfill \quad(definition of correctness).

    \item For any model $\mathfrak{M}=(\mathfrak{F},V)=(W,R,V)$,\\
    \indent\quad $\mathfrak{M}$ is a BIO-model
    
    $\Leftrightarrow$ for any $\alpha\in Form_{\to}$, $\Vert \alpha \Vert \subseteq \Box_{\mathfrak{F}}\Vert\alpha \Vert$ \hfill (definition of BIO-model, lemma \ref{formula.fixedpoint}, $R$ is transitive, lemma \ref{frame condition.fixpoint})
    
    $\Leftrightarrow$ for any $\alpha\in Form_{\to}$, $\Vert \alpha \Vert \subseteq \Vert \Box\alpha\Vert$  \hfill \quad(definition of $\Vert\cdot\Vert$, $\Vert \Box\alpha\Vert=\Box_{\mathfrak{F}}\Vert\alpha \Vert$)
    
    $\Leftrightarrow$ for any $\alpha\in Form_{\to}$, $ \alpha \,\vDash_{\mathfrak{M}} \Box\alpha$  \hfill \quad(definition of the satisfaction relation and $\,\vDash_{\mathfrak{M}}$)
    
    $\Leftrightarrow$ (Prop${}_{tr}$) is correct at $\mathfrak{M}$  \hfill \quad(definition of correctness).
   
\end{enumerate}
\end{proof}

Now, we prove soundness of the eight systems.
\begin{theorem}[Soundness Theorems]
 For each $\clubsuit\in\{\mathbf{K}^p, \mathbf{T}^p, \mathbf{B}^p, \mathbf{V}, \mathbf{B4}^p, \mathbf{I}, \mathbf{O}, \mathbf{C}\}$, for any $\Gamma\subseteq Form_\to$ and $\phi\in Form_\to$,
    \[\Gamma\,\vdash_\clubsuit\phi \Rightarrow\Gamma\,\vDash_\clubsuit\phi.\]
\end{theorem}

\begin{proof}
    Fix $\clubsuit\in\{\mathbf{K}^p, \mathbf{T}^p, \mathbf{B}^p, \mathbf{V}, \mathbf{B4}^p, \mathbf{I}, \mathbf{O}, \mathbf{C}\}$ and $\mathfrak{M}\in \mathcal{D}_\clubsuit$. We prove that $\,\vdash_\clubsuit\, \subseteq\, \,\vDash_{\mathfrak{M}}$. It suffices to show that $\,\vDash_{\mathfrak{M}}$ satisfies all the rules defining $\,\vdash_\clubsuit$. This readily follows from the previous five lemmas.
\end{proof}

\subsection{Completeness}

Let $\ \,\vdash\ \subseteq \wp(Form_\to) \times Form_\to$. Similar to Section 3.3, one can define the notions of being ``\(\,\vdash\)-consistent", ``\(\,\vdash\)-\(\phi\)-maximal consistent", ``\(\,\vdash\)-closed", and ``\(\,\vdash\)-prime". And similarly it can proved that being \(\,\vdash\)-\(\phi\)-maximal consistent implies being \(\,\vdash\)-consistent, \(\,\vdash\)-closed and \(\,\vdash\)-prime.

And analogously one can prove the Lindenbaum's lemma.
\begin{lemma}[Lindenbaum's lemma]
Let $\ \,\vdash\ \subseteq \wp(Form_\to) \times Form_\to$ satisfy \( (A), (Cut)\), \((Com), (Mon) \).
If \(\Gamma \nvdash\phi\), then there exists \(\Phi \subseteq \text{Form}_\to\), such that \(\Phi\) is \(\,\vdash\)-\(\phi\)-maximal consistent and \(\Gamma \subseteq \Phi\) and $\phi\notin\Phi$.
\end{lemma}

By virtue of the increase in expressive power of the language, we can have a uniform way of defining canonical models for various systems. (In Section 3.3, the definition of the canonical accessibility relation varies in different systems.)

\begin{definition}[Canonical model for $\,\vdash$]
Let $\ \,\vdash\ \subseteq \wp(Form_\to) \times Form_\to$. Define:
\begin{enumerate}[(i)]
    \item \(W^c_{\,\vdash} = \{ \Gamma \subseteq \text{Form}_\neg \mid \Gamma\) is \(\vdash\)-closed,  $\vdash$-consistent and  $\,\vdash$-prime\(\}\).
    \item $R^c_{ \,\vdash}=(W^c_{ \,\vdash})^2\cap R^c$, where \(R^c = \{ (\Gamma,\Delta)\in \wp(Form)^2 \mid \) for any \(\alpha,\beta\in Form_\to\),  if $\alpha\to\beta\in\Gamma$ and  $\alpha\in\Delta$, then \(\beta\in\Delta\}\).
    \item \(V^c_{ \,\vdash}\) is a function from  \(W^c_{\,\vdash}\) to \(\mathcal{P}(PL)\) such that \(V^c_{ \,\vdash}(\Gamma) = \{ p \in PL \mid p \in \Gamma \}\) for any $\Gamma\in W^c_{\,\vdash}$.
\end{enumerate}
Then define
\(
\mathcal{M}^c_{ \vdash} = ( \mathcal{F}^c_{ \,\vdash}, V^c_{ \,\vdash} )
\), where $\mathcal{F}^c_{ \,\vdash} = ( W^c_{ \,\vdash}, R^c_{ \,\vdash})$.
\end{definition}

Before we can prove the Existence Lemma, we need some technical preparations. 

\begin{definition}
For any $\Gamma,\Delta\in\wp(Form_\to)$, define  $\Gamma^\to(\Delta)=\{\psi\in Form\;|\; $there is $\alpha\in\Delta$ with $ \alpha\to\psi\in\Gamma\}$. 
\end{definition}

\begin{lemma}\label{Gamma(Delta,varphi)}
Let $\,\vdash\,\subseteq \wp(Form_\to) \times Form_\to$ satisfy all the basic rules and (Com). Let $\Gamma,\Delta\subseteq Form_\to$ be $\,\vdash$-closed. Then $\Delta\subseteq\Gamma^\to(\Delta)$, $\Gamma R^c \Gamma^\to(\Delta)$ and $\Gamma^\to(\Delta)$ is $\,\vdash$-closed.
\end{lemma}

\begin{proof}
~
\begin{enumerate}[(1)]
    \item For any $\psi\in\Delta$, we have $\psi\vdash\psi$ by (A), then $\Gamma\vdash \psi\to\psi$ by (DT${}_0$) and (Mon). Then $ \psi\to\psi\in\Gamma$ since $\Gamma$ is $\,\vdash$-closed. So $\psi\in \Gamma^\to(\Delta)$.
   
    \item For any $\phi,\psi\in Form$, if $\phi\to\psi\in\Gamma$ and $\phi\in\Gamma^\to(\Delta)$, that is,  $\phi\to\psi\in\Gamma$ and $ \alpha\to\phi\in\Gamma$ for some $\alpha\in\Delta$, then by ($\to$tr) and (Mon) we have $\Gamma\vdash \alpha\to\psi$, so $\psi\in\Gamma^\to(\Delta)$.
    
    \item Suppose $\Gamma^\to(\Delta)\vdash\psi$. By (Com),(Mon) and ($\land$E), we may assume that there is $ n\in\mathbb{N}^*$ and $\psi_1 ,...,\psi_n \in \Gamma^\to(\Delta)$ such that $\psi_1 \land...\land\psi_n\vdash\psi$. By definition of $\Gamma^\to(\Delta)$, there are $\alpha_1 ,..., \alpha_n\in\Delta$ such that $ \alpha_i \to\psi_i\in\Gamma$ for $i=1,...,n$. By ($\land$E) we can get $\alpha_1 \land...\land\alpha_n \vdash \alpha_i $ for each $i=1,...,n$. Then by (A),(DT${}_0$), ($\to$tr) and (Cut) we can obtain $\Gamma\vdash\alpha_1 \land...\land\alpha_n \to\psi_i$ for each $i=1,...,n$. Then $\Gamma\vdash\alpha_1 \land...\land\alpha_n \to\psi_1 \land...\land\psi_n$ by ($\to\land$). Together with $\psi_1 \land...\land\psi_n\vdash\psi$, using (DT${}_0$), ($\to$tr), (Cut) and closedness of $\Gamma$, we have  $\alpha_1 \land...\land\alpha_n \to\psi\in\Gamma$. Since $\alpha_1 ,..., \alpha_n\in\Delta$ and $\Delta$ is $\,\vdash$-closed, $\alpha_1 \land...\land\alpha_n \in \Delta$. So $\psi\in\Gamma^\to(\Delta)$. 
   
 \end{enumerate}
\end{proof}

\begin{lemma}[Existence lemma]
Let $\ \,\vdash\ \subseteq \wp(Form_\to) \times Form_\to$ satisfy all the basic rules. For any $\,\vdash$-closed $\Gamma \subseteq Form_\to$ and $\alpha,\beta \in Form_\to$, if $ \alpha\to\beta \notin \Gamma$, then there exists $ \Phi \in W^c_{ \,\vdash}$, such that $\Gamma R^c \Phi$, $\alpha\in \Phi$ and $\beta\notin \Phi$.
\end{lemma}

\begin{proof}
Let $\Omega = \{ \Delta \mid \Delta$ is $\vdash$-consistent and $\vdash$-closed, $\Gamma R^c \Delta$, $\alpha\in \Delta$ and $\beta\notin\Delta\}$.

Consider $\Gamma^\to(\text{Th}_{\,\vdash}(\{\alpha\}))$. Since $\text{Th}_{\,\vdash}(\{\alpha\})$ can be easily proved to be $\,\vdash$-closed, by the previous lemma, $\Gamma^\to(\text{Th}_{\,\vdash}(\{\alpha\}))$ is $\,\vdash$-closed, $\Gamma R^c \Gamma^\to(\text{Th}_{\,\vdash}(\{\alpha\}))$ and $\alpha\in \Gamma^\to(\text{Th}_{\,\vdash}(\{\alpha\}))$. Moreover, since $\alpha\to\beta \notin \Gamma$, we can show that $\beta\notin\Gamma^\to(\text{Th}_{\,\vdash}(\{\alpha\}))$, hence $\Gamma^\to(\text{Th}_{\,\vdash}(\{\alpha\}))$ is $\,\vdash$-consistent by ($\bot$). So $\Gamma^\to(\text{Th}_{\,\vdash}(\{\alpha\})) \in \Omega$ and $\Omega$ is not empty.

It is easy to show (using (Com), (Mon)) that for any directed family $\mathcal{H} \subseteq \Omega$, $\bigcup \mathcal{H} \in \Omega$.

By Zorn's lemma, $\Omega$ has a maximal element, denoted as $\Phi$. We prove that $\Phi$ is $\,\vdash$-prime.

Let $2 \leq n \in \mathbb{N}^*$, and $\phi, \alpha_1, \ldots, \alpha_n\in Form_\to$. Assume that $\phi \in \Phi$ and $\phi \vdash (\alpha_1, \ldots, \alpha_n)$.
 Suppose towards contradiction that none of $\alpha_1, \ldots, \alpha_n$ belongs to $\Phi$.

Consider $\Gamma^\to(\text{Th}_{\,\vdash}(\Phi \cup \{\alpha_i\}))$. Note that $\Gamma^\to(\text{Th}_{\,\vdash}(\Phi \cup \{\alpha_i\}))$ is $\,\vdash$-closed, $\Gamma R^c \Gamma^\to(\text{Th}_{\,\vdash}(\Phi \cup \{\alpha_i\}))$ and $\Phi\subset \Gamma^\to(\text{Th}_{\,\vdash}(\Phi \cup \{\alpha_i\}))$, by the previous lemma. And note that $\bot \in \Gamma^\to(\text{Th}_{\,\vdash}(\Phi \cup \{\alpha_i\}))$ implies $\beta\in\Gamma^\to(\text{Th}_{\,\vdash}(\Phi \cup \{\alpha_i\}))$. By the maximality of $\Phi$, we have $\beta\in\Gamma^\to(\text{Th}_{\,\vdash}(\Phi \cup \{\alpha_i\}))$. Namely, for any $i \in \{1, \ldots, n\}$, there exists $\chi_i\in\Gamma$ such that $\Phi \cup \{\alpha_i\}\,\vdash\chi_i$ and $\chi_i\to\beta\in\Gamma$.

Using (Com), (Mon), ($\land$E), (Cut) and closedness of $\Phi$, we can infer that there exists $\psi \in \Phi$, for any $i \in \{1, \ldots, n\}$,  $\psi \land \alpha_i\,\vdash\chi_i$.

Using (DT${}_0$), we have $\,\vdash\psi \land \alpha_i\to\chi_i$ for any $i \in \{1, \ldots, n\}$. Since $\chi_i\to\beta\in\Gamma$ and $\Gamma$ is $\,\vdash$-closed, using (A),($\to$tr) and (Cut), we get $\psi \land \alpha_i\to\beta \in \Gamma$ for any $i \in \{1, \ldots, n\}$.

Moreover, since $\phi \vdash (\alpha_1, \ldots, \alpha_n)$, we have $(\psi \land \alpha_1\to\beta) \land \ldots \land (\psi \land \alpha_n\to\beta) \vdash \psi \land \phi\to\beta$.

Thus, using (A), ($\land$I), (Cut) and closedness of $\Gamma$, we have $\psi \land \phi\to\beta \in \Gamma$.

On the other hand, since $\psi, \phi \in \Phi$ and $\beta\notin\Phi$, using ($\land$I) and closedness of $\Phi$, we have $\psi \land \phi \in \Phi$ and $\beta\notin\Phi$, contradicting $\Gamma R^c \Phi$.
\end{proof}

\begin{lemma}[Truth Lemma]
Let $\ \,\vdash\ \subseteq \wp(Form_\to) \times Form_\to$ satisfy all the basic rules. For any $\Gamma \in W^c_{ \,\vdash}$ and $\alpha \in Form_\to$, 
\[\mathfrak{M}^c_{ \,\vdash},\Gamma\,\vDash\alpha \iff \alpha\in\Gamma.\]
\end{lemma}

\begin{proof}
 Use induction on formulas to show that for any $\Gamma \in W^c_{ \,\vdash}$, the target equivalence holds. The case for propositional letters follows from the definition of $V^c_{\,\vdash}$. For any $\Gamma\in W^c_{\,\vdash}$, $\bot\notin\Gamma$, since $\Gamma$ is ${\,\vdash}$-consistent and $\,\vdash$ satisfies (A). The case for $\land$ is proved using ($\land$I), ($\land$E) and the closedness of each $\Gamma\in W^c_{\,\vdash}$. The case for $\to$ follows from the previous lemma and the fact that \(\alpha\to\beta\in\Gamma\) and $\alpha\in\Delta$ implies \(\beta\in\Delta\) for any $(\Gamma,\Delta)\in R^c_{\,\vdash}$ and $\alpha,\beta\in Form_\to$.
\end{proof}

The next two lemmas establish connections between the non-basic rules and the extra properties of the canonical model.
\begin{lemma}\label{property of CM 1}
Let $\ \,\vdash\ \subseteq \wp(Form_\to) \times Form_\to$ satisfy the basic rules.
\begin{enumerate}[(1)]
    \item If $\,\vdash$ satisfies (Refl), then $R^c_\vdash$ is reflexive.
    \item If $\,\vdash$ satisfies (Prop${}_{tr}$), then $R^c_\vdash$ is transitive.
    \item If $\,\vdash$ satisfies (Sym), then $R^c_\vdash$ is symmetric.
\end{enumerate}
\end{lemma}

\begin{proof}
~
\begin{enumerate}[(1)]
    \item Let $\Gamma\in W^c_\vdash$ be arbitrary. We want to show $\Gamma R^c_\vdash \Gamma$. For any formula $\varphi,\psi$, if $\varphi\to\psi\in\Gamma$ and $\varphi\in\Gamma$, then by (A),($\land$I),(Refl) and (Cut), we have $\Gamma\vdash\psi$, so $\psi\in\Gamma$ by the closedness of $\Gamma$. 
    
    \item Let $\Gamma,\Delta,\Phi \in W^c_\vdash$, $\Gamma R^c_\vdash \Delta$ and $\Delta R^c_\vdash\Phi $. We want to show $\Gamma R^c_\vdash \Phi$. For any formula $\varphi,\psi$, if $\varphi\to\psi\in\Gamma$, then $T\to (\varphi\to\psi)\in\Gamma$ by (Prop${}_{tr}$), (Mon) and closedness of $\Gamma$, then $\varphi\to\psi\in\Delta$ since $\Gamma R^c_\vdash \Delta$ and $T\in\Delta$ (by (DT${}_0$) and closedness of $\Delta$), then $\varphi\in\Phi$ implies $\psi\in\Phi$ since $\Delta R^c_\vdash\Phi $.
    
    \item Now, let $\Gamma,\Delta \in W^c_\vdash$ and $\Gamma R^c_\vdash \Delta$. We want to show $\Delta R^c_\vdash \Gamma$. Let $\varphi,\psi\in Form_\to$, $\varphi\to\psi\in\Delta$ and $\varphi\in\Gamma$. Suppose towards contradiction that $\psi\notin\Gamma$. Since $\,\vdash$ satisfies (Sym) and $\Gamma$ is $\,\vdash$-prime, we have $(\varphi\to\psi)\to\bot\in\Gamma$. Then $\bot\in\Delta$ since $\Gamma R^c_\vdash \Delta$ and $\varphi\to\psi\in\Delta$, which contradicts the consistency of $\Delta$. 
\end{enumerate}
\end{proof}

\begin{lemma}\label{property of CM 2}
Let $\ \,\vdash\ \subseteq \wp(Form_\to) \times Form_\to$ satisfy the basic rules.
\begin{enumerate}[(1)]
    \item If $\,\vdash$ satisfies (Prop${}_{\mathbf{K}^p-}$), then $\mathfrak{M}^c_\vdash\in \mathcal{D}_{\mathbf{K}^p-}$.
    \item If $\,\vdash$ satisfies (Refl), then $\mathfrak{M}^c_\vdash\in \mathcal{D}_{\mathbf{T}^p-}$.
    \item If $\,\vdash$ satisfies (Prop${}_{tr}$), then $\mathfrak{M}^c_\vdash\in \mathcal{D}_{\mathbf{V}}$.
    \item If $\,\vdash$ satisfies (Sym) and (Prop${}_{sy}$), then $\mathfrak{M}^c_\vdash\in \mathcal{D}_{\mathbf{B}^p}$.
    \item If $\,\vdash$ satisfies (Refl), (Sym) and (Prop${}_{sy}$), then $\mathfrak{M}^c_\vdash\in \mathcal{D}_{\mathbf{O}}$.
    \item If $\,\vdash$ satisfies (Sym) and (Prop${}_{tr}$), then $\mathfrak{M}^c_\vdash\in \mathcal{D}_{\mathbf{B4}^p}$.
    \item If $\,\vdash$ satisfies (Refl) and (Prop${}_{tr}$), then $\mathfrak{M}^c_\vdash\in \mathcal{D}_{\mathbf{I}}$.
    \item If $\,\vdash$ satisfies (Refl), (Sym) and (Prop${}_{tr}$), then $\mathfrak{M}^c_\vdash\in \mathcal{D}_{\mathbf{C}}$.
\end{enumerate}
\end{lemma}

\begin{proof}
~
\begin{enumerate}[(1)]
    \item \quad ~$\mathfrak{M}^c_\vdash\in \mathcal{D}_{\mathbf{K}^p-}$ \\
    $\Leftrightarrow$ (Prop${}_{\mathbf{K}^p-}$) is correct at $\mathfrak{M}^c_\vdash$ \hfill (lemma \ref{definability})

    $\Leftrightarrow$  for any $\alpha\in Form_\to$, $ \alpha  \,\vDash_{\mathfrak{M}^c_\vdash}  \neg T\to\alpha$ \hfill (definition of correctness)
    
    $\Leftrightarrow$ for any $\Gamma\in W^c_\vdash$, for any $\alpha\in Form_\to$, if $\alpha\in\Gamma$, then $ \neg T\to\alpha \in\Gamma$ \hfill (definition of $\,\vDash_{\mathfrak{M}^c_\vdash}$, Truth Lemma)
    
    $\Leftarrow$ $\,\vdash$ satisfies (Prop${}_{\mathbf{K}^p-}$)\hfill ((Mon), every $\Gamma\in W^c_\vdash$ is $\,\vdash$-closed)
    
    \item Since $\,\vdash$ satisfies (Refl), $R^c_\vdash$ is reflexive by lemma \ref{property of CM 1}. So, $\mathfrak{M}^c_\vdash\in \mathcal{D}_{\mathbf{T}^p-}$ by the definition of $\mathcal{D}_{\mathbf{T}^p-}$. 
    
    \item Since $\,\vdash$ satisfies (Prop${}_{tr}$), $R^c_\vdash$ is transitive by lemma \ref{property of CM 1}. So, we have:
    
    \indent\quad ~$\mathfrak{M}^c_\vdash\in \mathcal{D}_{\mathbf{V}}$ 
    
    $\Leftrightarrow$ (Prop${}_{tr}$) is correct at  $\mathfrak{M}^c_\vdash$ \hfill ($R^c_\vdash$ is transitive, lemma \ref{definability})

    $\Leftrightarrow$  for any $\alpha\in Form_\to$, $ \alpha  \,\vDash_{\mathfrak{M}^c_\vdash}   T\to\alpha$ \hfill (definition of correctness)
    
    $\Leftrightarrow$ for any $\Gamma\in W^c_\vdash$, for any $\alpha\in Form_\to$, if $\alpha\in\Gamma$, then $ T\to\alpha \in\Gamma$ \hfill (definition of $\,\vDash_{\mathfrak{M}^c_\vdash}$, Truth Lemma)
    
    $\Leftarrow$ $\,\vdash$ satisfies (Prop${}_{tr}$)\hfill ((Mon), every $\Gamma\in W^c_\vdash$ is $\,\vdash$-closed)

    \item Since $\,\vdash$ satisfies (Sym),  $R^c_\vdash$ is symmetric by lemma \ref{property of CM 1}. So, we have:
    
    \indent\quad ~$\mathfrak{M}^c_\vdash\in \mathcal{D}_{\mathbf{B}^p}$ 
    
    $\Leftrightarrow$ (Prop${}_{sy}$) is correct at  $\mathfrak{M}^c_\vdash$ \hfill ($R^c_\vdash$ is symmetric, lemma \ref{definability})

    $\Leftrightarrow$  for any $\alpha\in Form_\to$, $ \alpha  \,\vDash_{\mathfrak{M}^c_\vdash}   \neg\neg\alpha\to\alpha$ \hfill (definition of correctness)
    
    $\Leftrightarrow$ for any $\Gamma\in W^c_\vdash$, for any $\alpha\in Form_\to$, if $\alpha\in\Gamma$, then $ \neg\neg\alpha\to\alpha \in\Gamma$ \hfill (definition of $\,\vDash_{\mathfrak{M}^c_\vdash}$, Truth Lemma)
    
    $\Leftarrow$ $\,\vdash$ satisfies (Prop${}_{sy}$)\hfill ((Mon), every $\Gamma\in W^c_\vdash$ is $\,\vdash$-closed)
    
    \item[(5)$\sim$(8)] Note that: \\
    \indent\quad $\mathcal{D}_{\mathbf{O}}=\{\mathfrak{M}\in\mathcal{D}_{\mathbf{B}^p}\,|\, \mathfrak{M}$ is reflexive $\}$,\\
    \indent\quad $\mathcal{D}_{\mathbf{B4}^p}=\{\mathfrak{M}\in\mathcal{D}_{\mathbf{V}}\,|\, \mathfrak{M}$ is symmetric $\}$,\\
    \indent\quad $\mathcal{D}_{\mathbf{I}}=\{\mathfrak{M}\in\mathcal{D}_{\mathbf{V}}\,|\, \mathfrak{M}$ is reflexive $\}$,\\
    \indent\quad $\mathcal{D}_{\mathbf{C}}=\{\mathfrak{M}\in\mathcal{D}_{\mathbf{V}}\,|\, \mathfrak{M}$ is reflexive and symmetric $\}$.\\
    So (5)$\sim$(8) follows from (3), (4) and lemma \ref{property of CM 1}.
    
\end{enumerate}
\end{proof}

Now, we prove the completeness of the eight systems.
\begin{theorem}[Completeness Theorems]
  For each $\clubsuit\in\{\mathbf{K}^p,\mathbf{T}^p, \mathbf{B}^p, \mathbf{V}$,$\mathbf{B4}^p$, $\mathbf{I}$,$\mathbf{O}, \mathbf{C}\}$, for any $\Gamma\subseteq Form_\to$ and $\phi\in Form_\to$,
    \[\Gamma\,\vDash_\clubsuit\phi \Rightarrow\Gamma\,\vdash_\clubsuit\phi.\]
\end{theorem}
\begin{proof}
    Let $\clubsuit\in\{\mathbf{K}^p,\mathbf{T}^p, \mathbf{B}^p, \mathbf{V}, \mathbf{B4}^p, \mathbf{I}, \mathbf{O}, \mathbf{C}\}$. Let $\Gamma\subseteq Form_\to$ and $\phi\in Form_\to$.We prove the contrapositive of the target entailment. 
    
    Assume that $\Gamma\,\nvdash_\clubsuit\phi$. Then by Lindenbaum's lemma, there exists \(\Phi \in W^c_{ \,\vdash_\clubsuit}\) with \(\Gamma \subseteq \Phi\) and $\phi\notin\Phi$. By the Truth Lemma, $\mathfrak{M}^c_{ \,\vdash_\clubsuit},\Phi\,\vDash\Gamma$ and $\mathfrak{M}^c_{ \,\vdash_\clubsuit},\Phi\,\nvDash\alpha$. Moreover, according to lemma \ref{property of CM 2}, $\mathfrak{M}^c_{ \,\vdash_\clubsuit}\in \mathcal{D}_\clubsuit$ for $\clubsuit\in\{\mathbf{V}, \mathbf{B}^p, \mathbf{B4}^p, \mathbf{I}, \mathbf{O}, \mathbf{C}\}$, while $\mathfrak{M}^c_{ \,\vdash_{\clubsuit}}\in \mathcal{D}_{\clubsuit-}$ for $\clubsuit\in\{\mathbf{K}^p, \mathbf{T}^p\}$. So $\Gamma\,\nvDash_\clubsuit\phi$ if $\clubsuit\in\{\mathbf{V}, \mathbf{B}^p, \mathbf{B4}^p,\mathbf{I}, \mathbf{O}, \mathbf{C}\}$, and $\Gamma\,\nvDash_{\clubsuit-}\phi$ if $\clubsuit\in\{\mathbf{K}^p, \mathbf{T}^p\}$. Finally, by lemma \ref{vDash_K-=vDash_K, vDash_T-=vDash_T}, we have $\Gamma\,\nvDash_{\clubsuit}\phi$ if $\clubsuit\in\{\mathbf{K}^p, \mathbf{T}^p\}$.
\end{proof}

\section{Translation into modal logics}
First, let us review some notions in Kripke semantics for modal logic.

\begin{definition}
Let $Form_\Box$ denote the set of modal formulas, which are built from propositional letters, using connectives $\neg,\land$ and modal operator $\Box$. (Moreover, for modal formulas, let $\bot$, $\alpha\lor\beta$, $\alpha\to\beta$, $\Diamond\alpha$ be the abbreviations of $p_0\land\neg p_0$, $\neg(\neg\alpha\land\neg\beta)$, $\neg(\alpha\land\neg\beta)$, and $\neg\Box\neg\alpha$, respectively.)

The satisfaction relation $\Vdash$ for modal formulas is recursively defined as follows. For any pointed Kripke model $\mathfrak{M},w=(W,R,V,w)$, let 
\begin{itemize}
    \item $\mathfrak{M},w\Vdash p\iff w\in V(p)$; 
    \item $\mathfrak{M},w\Vdash \alpha\land\beta\iff\mathfrak{M},w\Vdash \varphi$ and $\mathfrak{M},w\Vdash \beta$;
    \item $\mathfrak{M},w\Vdash \Box\alpha\iff$  for each $v\in W$, if $wRv$ then $\mathfrak{M},v\Vdash \alpha$.
\end{itemize}

For any modal formula $\alpha$ and Kripke model $\mathfrak{M}$, define the truth set of $\alpha$ (as a modal formula) in $\mathfrak{M}$ to be $\Vert\alpha\Vert^\Vdash_\mathfrak{M}=\{w\in\mathfrak{M}\mid \mathfrak{M},w\Vdash\alpha\}$.

A modal formula $\alpha$ is valid in a model $\mathfrak{M}$, denoted by $\mathfrak{M}\Vdash\alpha$, if and only if $\mathfrak{M},w\Vdash\alpha$ for any $w\in\mathfrak{M}$.

For any class $\mathcal{D}$ of Kripke models, the modal semantic consequence relation with respect to $\mathcal{D}$ is defines as: for any $\Gamma\subseteq Form_\Box$ any $\phi\in Form_\Box$,
$$\Gamma \Vdash_\mathcal{D} \phi \iff \text{for any } \mathfrak{M} \in \mathcal{D}, \text{ for any } w \in \mathfrak{M}, (\mathfrak{M}, w \Vdash \Gamma \Rightarrow \mathfrak{M}, w \Vdash \phi).$$
\end{definition}

Next, we define the translation functions and prove their properties within a relatively general framework.

\begin{definition}
Let $\varphi\in Form_\Box$. Let $\varphi(\cdot)$ be the unary operator on $Form_\Box$, such that for each $\alpha\in Form_\Box$, $$\varphi(\alpha)=[\alpha/p_0]\varphi,$$ 
where $[\alpha/p_0]\varphi$ is the result of substituting $p_0$ for $\alpha$ in $\varphi$.

Define recursively as follows the translation function (induced by $\varphi$) from $Form$, the set of propositional formulas, to $Form_\Box$.
\begin{itemize}
    \item $p^{\varphi}=\varphi(p)$; 
    \item $\bot^{\varphi}=\bot$;
    \item $(\alpha\land\beta)^{\varphi}=\alpha^{\varphi}\land \beta^{\varphi}$;
    \item $(\neg\alpha)^{\varphi}=\Box\neg(\alpha^{\varphi})$
    \item $(\alpha\to\beta)^{\varphi}=\Box(\alpha^{\varphi}\to \beta^{\varphi})$.
\end{itemize}
\end{definition}

\begin{definition}
    For any class $\mathcal{D}$ of Kripke models, for any $\varphi\in Form_\Box$, define $\mathcal{D}_\varphi=\{\mathfrak{M}\in\mathcal{D} \mid$ for any $p\in PL,\ \mathfrak{M}\Vdash p\leftrightarrow\varphi(p)\}$.
\end{definition}

\begin{lemma}[Tranlation lemma]\label{translation lemma}
    Let $\mathcal{D}$ be a class of Kripke models, and $\varphi\in Form_\Box$.  If $\,\mathfrak{M}\Vdash \varphi(p)\leftrightarrow\varphi(\varphi(p))$ for any $\mathfrak{M}\in\mathcal{D}$ and $p \in PL$, then for any $\Gamma\subseteq Form$ and $\phi\in Form$, 
    \[ \Gamma \vDash_{\mathcal{D}_\varphi} \phi \iff \Gamma^\varphi \Vdash_\mathcal{D} \phi^\varphi, \]
    where $\Gamma^\varphi=\{\alpha^\varphi\mid\alpha\in\Gamma\}$.
\end{lemma}

\begin{proof}
   Assume $\mathfrak{M}\Vdash \varphi(p)\leftrightarrow\varphi(\varphi(p))$ for any $\mathfrak{M}\in\mathcal{D}$ and $p \in PL$. It is sufficient to provide $F: \mathcal{D}\to\mathcal{D}_\varphi$ and $G: \mathcal{D}_\varphi\to\mathcal{D}$ such that: 
   \begin{align*}
       w\in F(\mathfrak{M}) &\,\textit{ and }\, v\in G(\mathfrak{N});\\
       F(\mathfrak{M}),w\vDash\alpha &\iff \mathfrak{M},w\Vdash\alpha^\varphi; \\
       \mathfrak{N},v\vDash\alpha &\iff G(\mathfrak{N},v)\Vdash\alpha^\varphi.
   \end{align*}
for any $\alpha\in Form$, $\mathfrak{M}\in\mathcal{D}$, $\mathfrak{N}\in\mathcal{D}_\varphi$, $w\in\mathfrak{M}$ and $v\in\mathfrak{N}$.

Define $F$ as: \(F(\mathfrak{M})=(W,R,\{(p,\Vert\varphi(p)\Vert^\Vdash_\mathfrak{M})\mid p\in PL\})\), for any $\mathfrak{M}=(W,R,V)\in\mathcal{D}$. Then it is easy to show that for any $\mathfrak{M}\in\mathcal{D}$, $w\in\mathfrak{M}$ and $p\in PL$, we have: (1) $w\in F(\mathfrak{M})$; (2) $F(\mathfrak{M}),w\Vdash p\Leftrightarrow \mathfrak{M},w\Vdash \varphi(p)$; (3) $F(\mathfrak{M}),w\Vdash \varphi(p)\Leftrightarrow \mathfrak{M},w\Vdash \varphi(\varphi(p))$. So $ran(F)\subseteq\mathcal{D}_\varphi$, since $\mathfrak{M}\Vdash \varphi(p)\leftrightarrow\varphi(\varphi(p))$ for any $\mathfrak{M}\in\mathcal{D}$. By induction on $Form$ we can readily prove that $F(\mathfrak{M}),w\vDash\alpha \Leftrightarrow \mathfrak{M},w\Vdash\alpha^\varphi$ for any $\alpha\in Form$.

Define $G$ as: \(G(\mathfrak{N})=\mathfrak{N}\) for any $\mathfrak{N}\in\mathcal{D}_\varphi$. Then $ran(G)\subseteq\mathcal{D}$, since $\mathcal{D}_\varphi\subseteq\mathcal{D}$. Clearly $v\in G(\mathfrak{N})$ for any $v\in\mathfrak{N}$. Since $\mathfrak{N}\in\mathcal{D}_\varphi$, $\mathfrak{N}\Vdash p\leftrightarrow\varphi(p)$. Then by induction we can easily check that $\mathfrak{N},v\vDash\alpha \Leftrightarrow G(\mathfrak{N}),v\Vdash\alpha^\varphi$ for any $\alpha\in Form$.
\end{proof}

\begin{definition}
    Define $\varphi_\spadesuit\in Form$ for each $\spadesuit\in\{\mathbf{K},\mathbf{T}, \mathbf{B}, \mathbf{4}, \mathbf{B4}, \mathbf{S4}, \mathbf{TB}, \mathbf{S5}\}$ as follows.
    \begin{itemize}
        \item $\varphi_{\mathbf{K}}=p_0\land\Box(\Box\bot\to p_0)$;
        \item $\varphi_{\mathbf{T}}=p_0$;
        \item $\varphi_{\mathbf{B}}=p_0\land\Box(\Box\Diamond p_0\to p_0)$;
        \item $\varphi_{\mathbf{4}}=\varphi_{\mathbf{B4}}=p_0\land\Box p_0$;
        \item $\varphi_{\mathbf{TB}}=\Box\Diamond p_0$;
        \item $\varphi_{\mathbf{S4}}=\varphi_{\mathbf{S5}}=\Box p_0$.
    \end{itemize}
\end{definition}

Recall the classes $\mathcal{D}_{\spadesuit}$ and $\mathcal{D}_{\spadesuit^p}$ defined in Section 2. We establish the following connection through $\varphi_\spadesuit$:

\begin{lemma}\label{D_varphi connection}
    For any $\spadesuit\in\{\mathbf{B}, \mathbf{4}, \mathbf{B4}, \mathbf{S4}, \mathbf{TB}, \mathbf{S5}\}$,  
    $$\mathcal{D}_{\spadesuit^p}=(\mathcal{D}_{\spadesuit})_{\varphi_\spadesuit}.$$
    Besides, $\mathcal{D}_{\spadesuit^p-}=(\mathcal{D}_{\spadesuit})_{\varphi_\spadesuit}$ for $\spadesuit\in\{\mathbf{K},\mathbf{T}\}$.
\end{lemma}
\begin{proof}
    By definition, we need to prove that $\mathcal{D}_{\spadesuit^p}$(or $\mathcal{D}_{\spadesuit^p-}$, when $\spadesuit\in\{\mathbf{K},\mathbf{T}\}$) equals $\{\mathfrak{M}\in\mathcal{D}_{\spadesuit} \mid$ for any $p\in PL,\ \mathfrak{M}\Vdash p\leftrightarrow\varphi_{\spadesuit}(p)\}$. This is not difficult with the help of lemma \ref{frame condition.fixpoint}.
\end{proof}

Moreover, we can check that $\varphi_\spadesuit$ satisfies the requirement mentioned in the premise of lemma \ref{translation lemma}, 

\begin{lemma}\label{check requirement on varphi}
    For any $\spadesuit\in\{\mathbf{K}, \mathbf{T},\mathbf{B}, \mathbf{4}, \mathbf{B4}, \mathbf{S4}, \mathbf{TB}, \mathbf{S5}\}$, for any $\mathfrak{M}\in\mathcal{D}_\spadesuit$ and $p \in PL$,
    $$\mathfrak{M}\Vdash \varphi_\spadesuit(p)\leftrightarrow\varphi_\spadesuit(\varphi_\spadesuit(p)).$$
\end{lemma}
\begin{proof}
    Note that the target is equivalent to $\Vert \varphi_\spadesuit(p)\Vert^\Vdash_{\mathfrak{M}}=\Vert\varphi_\spadesuit(\varphi_\spadesuit(p))\Vert^\Vdash_{\mathfrak{M}}$. We only prove the cases for $\spadesuit\in\{\mathbf{K}, \mathbf{B}, \mathbf{4}\}$, for the rest are simple. We also omit the superscript and subscript of $\Vert\cdot\Vert$ for simplicity. Let $\mathfrak{M}=(\mathfrak{F},V)=(W,R,V)\in\mathcal{D}_\spadesuit$ and $p\in PL$.
    \begin{enumerate}[(1)]
        \item $\spadesuit=\mathbf{K}$. Note that $\Vert\varphi_\mathbf{K}(\varphi_\mathbf{K}(p))\Vert=\Vert \varphi_\mathbf{K}(p)\Vert\cap\Box_\mathfrak{F}\Vert \Box\bot\to\varphi_\mathbf{K}(p)\Vert$, so it suffices to show that $\Vert \varphi_\mathbf{K}(p)\Vert\subseteq\Box_\mathfrak{F}\Vert \Box\bot\to\varphi_\mathbf{K}(p)\Vert$. This is true since $\Vert \Box\bot\to\varphi_\mathbf{K}(p)\Vert\subseteq\Vert \Box\bot\to\Box(\Box\bot\to p)\Vert=W$.
        
        \item $\spadesuit=\mathbf{B}$. Note that $\Vert\varphi_\mathbf{B}(\varphi_\mathbf{B}(p))\Vert=\Vert \varphi_\mathbf{B}(p)\Vert\cap\Box_\mathfrak{F}\Vert \Box\Diamond\varphi_\mathbf{B}(p)\to\varphi_\mathbf{B}(p)\Vert$. So it suffices to show that $\Vert \varphi_\mathbf{B}(p)\Vert\subseteq\Box_\mathfrak{F}\Vert \Box\Diamond\varphi_\mathbf{B}(p)\to\varphi_\mathbf{B}(p)\Vert$, which is equivalent to $R\Vert \varphi_\mathbf{B}(p)\Vert\cap \Vert\Box\Diamond\varphi_\mathbf{B}(p)\Vert \subseteq \Vert \varphi_\mathbf{B}(p)\Vert$. We split it into two parts.
        \begin{enumerate}[(a)]
            \item $R\Vert \varphi_\mathbf{B}(p)\Vert\cap \Vert\Box\Diamond\varphi_\mathbf{B}(p)\Vert \subseteq \Vert p\Vert$ : \\
            Note that $R\Vert \varphi_\mathbf{B}(p)\Vert\subseteq R\Box_\mathfrak{F}\Vert \Box\Diamond p\to p\Vert\subseteq\Vert \Box\Diamond p\to p\Vert$. So $R\Vert \varphi_\mathbf{B}(p)\Vert\cap \Vert\Box\Diamond\varphi_\mathbf{B}(p)\Vert \subseteq \Vert \Box\Diamond p\to p\Vert\cap \Vert\Box\Diamond p\Vert \subseteq \Vert p\Vert$.

            \item $R\Vert \varphi_\mathbf{B}(p)\Vert\cap \Vert\Box\Diamond\varphi_\mathbf{B}(p)\Vert \subseteq \Box_\mathfrak{F}\Vert \Box\Diamond p\to p\Vert$ :\\
            Since $R$ is symmetric, and since $R\Vert \varphi_\mathbf{B}(p)\Vert\subseteq\Vert \Box\Diamond p\to p\Vert$ as proved in (a), $\Vert\Box\Diamond\varphi_\mathbf{B}(p)\Vert = \Box_\mathfrak{F}R\Vert\varphi_\mathbf{B}(p)\Vert\subseteq \Box_\mathfrak{F}\Vert \Box\Diamond p\to p\Vert$.    
        \end{enumerate}

        \item $\spadesuit=\mathbf{4}$. Note that $\Vert\varphi_\mathbf{4}(\varphi_\mathbf{4}(p))\Vert=\Vert \varphi_\mathbf{4}(p)\Vert\cap\Box_\mathfrak{F}\Vert \varphi_\mathbf{4}(p)\Vert=\Vert \varphi_\mathbf{4}(p)\Vert\cap\Box_\mathfrak{F}\Vert p\Vert\cap\Box_\mathfrak{F}\Box_\mathfrak{F}\Vert p\Vert$, so it suffices to show that $\Vert \varphi_\mathbf{4}(p)\Vert\subseteq\Box_\mathfrak{F}\Vert p\Vert\cap\Box_\mathfrak{F}\Box_\mathfrak{F}\Vert p\Vert$. This is easily verified using the transitivity of $R$.
    \end{enumerate}
\end{proof}

Combining lemma \ref{translation lemma}, lemma \ref{D_varphi connection}, lemma \ref{check requirement on varphi} and lemma \ref{vDash_K-=vDash_K, vDash_T-=vDash_T}, we have 
\begin{theorem}
    For any $\spadesuit\in\{\mathbf{K}, \mathbf{T},\mathbf{B}, \mathbf{4}, \mathbf{B4}, \mathbf{S4}, \mathbf{TB}, \mathbf{S5}\}$, for any $\Gamma\subseteq Form$ and $\phi\in Form$, 
    \[ \Gamma \vDash_{\mathcal{D}_{\spadesuit^p}} \phi \iff \Gamma^{\varphi_\spadesuit} \Vdash_{\mathcal{D}_\spadesuit} \phi^{\varphi_\spadesuit}.\]
\end{theorem}

\end{document}